% processed by citealice (July 21, 2000) on Fri Dec 15 15:34:57 CST 2000

% To: Saharon Shelah <shelah@math.huji.ac.il>,        Andrzej Roslanowski <roslanow@e-math.ams.org>
% Subject: sh703
% Date: Fri, 15 Dec 2000 15:29:26 -0500 (EST)
% From: Alice Leonhardt <leonhard@math.rutgers.edu>
% Mime-Version: 1.0
% Content-Description: revisions
% X-sliced-and-diced-by: 'savemail' 0.4, Jul 2000

% Latest revision - 00/Dec/18
% First typed - 99/Jan/21
% Previous version - 00/Nov/4

\ifx\shlhetal\undefinedcontrolsequence\let\shlhetal\relax\fi

\input amstex
% % \input mathdefs
%  *** start including mathdefs.tex *** 
\expandafter\ifx\csname mathdefs.tex\endcsname\relax
  \expandafter\gdef\csname mathdefs.tex\endcsname{}
\else \message{Hey!  Apparently you were trying to
  \string\input{mathdefs.tex} twice.   This does not make sense.} 
\errmessage{Please edit your file (probably \jobname.tex) and remove
any duplicate ``\string\input'' lines}\endinput\fi

%mathdefs.tex v1.3.2

%%% Changes from v1.0: footnote macros, warning for duplicated tags,
%%%   control sequences \( and \verbatimtags.
%%% From v1.2: \pretags, redefinition of \( using \ifinner, multi-part
%%%   equation numbering, control sequences \[, \references, and
%%%   \resetbracket. 
%%% From v1.3: \rm in \lastpart; write root of multi-part tag to .tgs 

%See file texdefs.doc for documentation.

\catcode`\X=12\catcode`\@=11

%Minor control sequences:
\def\n@wcount{\alloc@0\count\countdef\insc@unt}
\def\n@wwrite{\alloc@7\write\chardef\sixt@@n}
\def\n@wread{\alloc@6\read\chardef\sixt@@n}
\def\r@s@t{\relax}\def\v@idline{\par}\def\@mputate#1/{#1}
\def\l@c@l#1X{\firstpart.#1}\def\gl@b@l#1X{#1}\def\t@d@l#1X{{}}

%Creation of tag families and output of assignments and citations:
\def\crossrefs#1{\ifx\all#1\let\tr@ce=\all\else\def\tr@ce{#1,}\fi
   \n@wwrite\cit@tionsout\openout\cit@tionsout=\jobname.cit 
   \write\cit@tionsout{\tr@ce}\expandafter\setfl@gs\tr@ce,}
\def\setfl@gs#1,{\def\@{#1}\ifx\@\empty\let\next=\relax
   \else\let\next=\setfl@gs\expandafter\xdef
   \csname#1tr@cetrue\endcsname{}\fi\next}
\def\m@ketag#1#2{\expandafter\n@wcount\csname#2tagno\endcsname
     \csname#2tagno\endcsname=0\let\tail=\all\xdef\all{\tail#2,}
   \ifx#1\l@c@l\let\tail=\r@s@t\xdef\r@s@t{\csname#2tagno\endcsname=0\tail}\fi
   \expandafter\gdef\csname#2cite\endcsname##1{\expandafter
     \ifx\csname#2tag##1\endcsname\relax?\else\csname#2tag##1\endcsname\fi
     \expandafter\ifx\csname#2tr@cetrue\endcsname\relax\else
     \write\cit@tionsout{#2tag ##1 cited on page \folio.}\fi}
   \expandafter\gdef\csname#2page\endcsname##1{\expandafter
     \ifx\csname#2page##1\endcsname\relax?\else\csname#2page##1\endcsname\fi
     \expandafter\ifx\csname#2tr@cetrue\endcsname\relax\else
     \write\cit@tionsout{#2tag ##1 cited on page \folio.}\fi}
   \expandafter\gdef\csname#2tag\endcsname##1{\expandafter
      \ifx\csname#2check##1\endcsname\relax
      \expandafter\xdef\csname#2check##1\endcsname{}%
      \else\immediate\write16{Warning: #2tag ##1 used more than once.}\fi
      \multit@g{#1}{#2}##1/X%
      \write\t@gsout{#2tag ##1 assigned number \csname#2tag##1\endcsname\space
      on page \number\count0.}%
   \csname#2tag##1\endcsname}}

\def\multit@g#1#2#3/#4X{\def\t@mp{#4}\ifx\t@mp\empty%
      \global\advance\csname#2tagno\endcsname by 1 
      \expandafter\xdef\csname#2tag#3\endcsname
      {#1\number\csname#2tagno\endcsnameX}%
   \else\expandafter\ifx\csname#2last#3\endcsname\relax
      \expandafter\n@wcount\csname#2last#3\endcsname
      \global\advance\csname#2tagno\endcsname by 1 
      \expandafter\xdef\csname#2tag#3\endcsname
      {#1\number\csname#2tagno\endcsnameX}
      \write\t@gsout{#2tag #3 assigned number \csname#2tag#3\endcsname\space
      on page \number\count0.}\fi
   \global\advance\csname#2last#3\endcsname by 1
   \def\t@mp{\expandafter\xdef\csname#2tag#3/}%
   \expandafter\t@mp\@mputate#4\endcsname
   {\csname#2tag#3\endcsname\lastpart{\csname#2last#3\endcsname}}\fi}
\def\t@gs#1{\def\all{}\m@ketag#1e\m@ketag#1s\m@ketag\t@d@l p
\let\realscite\scite
\let\realstag\stag
   \m@ketag\gl@b@l r \n@wread\t@gsin
   \openin\t@gsin=\jobname.tgs \re@der \closein\t@gsin
   \n@wwrite\t@gsout\openout\t@gsout=\jobname.tgs }
\outer\def\localtags{\t@gs\l@c@l}
\outer\def\globaltags{\t@gs\gl@b@l}
\outer\def\newlocaltag#1{\m@ketag\l@c@l{#1}}
\outer\def\newglobaltag#1{\m@ketag\gl@b@l{#1}}

%Reading in tag information:
\newif\ifpr@ 
\def\m@kecs #1tag #2 assigned number #3 on page #4.%
   {\expandafter\gdef\csname#1tag#2\endcsname{#3}
   \expandafter\gdef\csname#1page#2\endcsname{#4}
   \ifpr@\expandafter\xdef\csname#1check#2\endcsname{}\fi}
\def\re@der{\ifeof\t@gsin\let\next=\relax\else
   \read\t@gsin to\t@gline\ifx\t@gline\v@idline\else
   \expandafter\m@kecs \t@gline\fi\let \next=\re@der\fi\next}
\def\pretags#1{\pr@true\pret@gs#1,,}
\def\pret@gs#1,{\def\@{#1}\ifx\@\empty\let\n@xtfile=\relax
   \else\let\n@xtfile=\pret@gs \openin\t@gsin=#1.tgs \message{#1} \re@der 
   \closein\t@gsin\fi \n@xtfile}

%Sections and subsections; local numbering:
\newcount\sectno\sectno=0\newcount\subsectno\subsectno=0
\newif\ifultr@local \def\ultralocal{\ultr@localtrue}
\def\firstpart{\number\sectno}
\def\lastpart#1{\ifcase#1 \or a\or b\or c\or d\or e\or f\or g\or h\or 
   i\or k\or l\or m\or n\or o\or p\or q\or r\or s\or t\or u\or v\or w\or 
   x\or y\or z \fi}

\def\resetall{\global\advance\sectno by 1\subsectno=0
   \gdef\firstpart{\number\sectno}\r@s@t}
\def\resetsub{\global\advance\subsectno by 1
   \gdef\firstpart{\number\sectno.\number\subsectno}\r@s@t}
\def\newsection#1\par{\resetall\vskip0pt plus.3\vsize\penalty-250
   \vskip0pt plus-.3\vsize\bigskip\bigskip
   \message{#1}\leftline{\bf#1}\nobreak\bigskip}
\def\subsection#1\par{\ifultr@local\resetsub\fi
   \vskip0pt plus.2\vsize\penalty-250\vskip0pt plus-.2\vsize
   \bigskip\smallskip\message{#1}\leftline{\bf#1}\nobreak\medskip}

%jj tags:
% On Andrzej's request:  we want to be able 
% to show tags as in noverbatim, with verbatim in the margin,
% and cites as in verbatim, with nonverbatim in the margin
% mg -- July 2000

\newdimen\marginshift

\newdimen\margindelta
\newdimen\marginmax
\newdimen\marginmin

\def\margininit{       
\marginmax=3 true cm                  % how much room, approximately
				      
\margindelta=0.1 true cm              % distance between entries
\marginmin=0.1true cm                 % where will leftmost entry be
\marginshift=\marginmin
}    % we cannot execute this right now, since 
     % there may be a \magnification coming later in the 
     % main file.   So we call \margininit at the end of 
     % alice2jlem

\def\t@gsjj#1,{\def\@{#1}\ifx\@\empty\let\next=\relax\else\let\next=\t@gsjj
   \def\@@{p}\ifx\@\@@\else
   \expandafter\gdef\csname#1cite\endcsname##1{\citejj{##1}}
   \expandafter\gdef\csname#1page\endcsname##1{?}
   \expandafter\gdef\csname#1tag\endcsname##1{\tagjj{##1}}\fi\fi\next}
\newif\ifshowstuffinmargin
\showstuffinmarginfalse
\def\jjtags{\showstuffinmargintrue
\ifx\all\relax\else\expandafter\t@gsjj\all,\fi}

\def\tagjj#1{\realstag{#1}\mginpar{\zeigen{#1}}}
\def\citejj#1{\zeigen{#1}\mginpar{\rechnen{#1}}}

\def\rechnen#1{\expandafter\ifx\csname stag#1\endcsname\relax ??\else
                           \csname stag#1\endcsname\fi}

\newdimen\theight

\def\marginfont{\sevenrm}

\def\trymarginbox#1{\setbox0=\hbox{\marginfont\hskip\marginshift #1}%
		\global\marginshift\wd0 
		\global\advance\marginshift\margindelta}

\def \mginpar#1{%
\ifvmode\setbox0\hbox to \hsize{\hfill\rlap{\marginfont\quad#1}}%
\ht0 0cm
\dp0 0cm
\box0\vskip-\baselineskip
\else 
             \vadjust{\trymarginbox{#1}%
		\ifdim\marginshift>\marginmax \global\marginshift\marginmin
			\trymarginbox{#1}%
                \fi
             \theight=\ht0
             \advance\theight by \dp0    \advance\theight by \lineskip
             \kern -\theight \vbox to \theight{\rightline{\rlap{\box0}}%
\vss}}\fi}

% \def\mginpar#1{mg-#1-mg }

%Verbatim tags:
\def\t@gsoff#1,{\def\@{#1}\ifx\@\empty\let\next=\relax\else\let\next=\t@gsoff
   \def\@@{p}\ifx\@\@@\else
   \expandafter\gdef\csname#1cite\endcsname##1{\zeigen{##1}}
   \expandafter\gdef\csname#1page\endcsname##1{?}
   \expandafter\gdef\csname#1tag\endcsname##1{\zeigen{##1}}\fi\fi\next}
\def\verbatimtags{\showstuffinmarginfalse
\ifx\all\relax\else\expandafter\t@gsoff\all,\fi}
\def\zeigen#1{\hbox{$\langle$}#1\hbox{$\rangle$}}
\def\margincite#1{\ifshowstuffinmargin\mginpar{\rechnen{#1}}\fi}

%Equation numbering:
\def\(#1){\edef\dot@g{\ifmmode\ifinner(\hbox{\noexpand\etag{#1}})
   \else\noexpand\eqno(\hbox{\noexpand\etag{#1}})\fi
   \else(\noexpand\ecite{#1})\fi}\dot@g}

%Reference numbering:
\newif\ifbr@ck
\def\eat#1{}
\def\[#1]{\br@cktrue[\br@cket#1'X]}
\def\br@cket#1'#2X{\def\temp{#2}\ifx\temp\empty\let\next\eat
   \else\let\next\br@cket\fi
   \ifbr@ck\br@ckfalse\br@ck@t#1,X\else\br@cktrue#1\fi\next#2X}
\def\br@ck@t#1,#2X{\def\temp{#2}\ifx\temp\empty\let\neext\eat
   \else\let\neext\br@ck@t\def\temp{,}\fi
   \def\teemp{#1}\ifx\teemp\empty\else\rcite{#1}\fi\temp\neext#2X}
\def\resetbr@cket{\gdef\[##1]{[\rtag{##1}]}}
\def\references{\resetbr@cket\newsection References\par}

%Footnotes:
\newtoks\symb@ls\newtoks\s@mb@ls\newtoks\p@gelist\n@wcount\ftn@mber
    \ftn@mber=1\newif\ifftn@mbers\ftn@mbersfalse\newif\ifbyp@ge\byp@gefalse
\def\defm@rk{\ifftn@mbers\n@mberm@rk\else\symb@lm@rk\fi}
\def\n@mberm@rk{\xdef\m@rk{{\the\ftn@mber}}%
    \global\advance\ftn@mber by 1 }
\def\rot@te#1{\let\temp=#1\global#1=\expandafter\r@t@te\the\temp,X}
\def\r@t@te#1,#2X{{#2#1}\xdef\m@rk{{#1}}}
\def\b@@st#1{{$^{#1}$}}\def\str@p#1{#1}
\def\symb@lm@rk{\ifbyp@ge\rot@te\p@gelist\ifnum\expandafter\str@p\m@rk=1 
    \s@mb@ls=\symb@ls\fi\write\f@nsout{\number\count0}\fi \rot@te\s@mb@ls}
\def\byp@ge{\byp@getrue\n@wwrite\f@nsin\openin\f@nsin=\jobname.fns 
    \n@wcount\currentp@ge\currentp@ge=0\p@gelist={0}
    \re@dfns\closein\f@nsin\rot@te\p@gelist
    \n@wread\f@nsout\openout\f@nsout=\jobname.fns }
\def\m@kelist#1X#2{{#1,#2}}
\def\re@dfns{\ifeof\f@nsin\let\next=\relax\else\read\f@nsin to \f@nline
    \ifx\f@nline\v@idline\else\let\t@mplist=\p@gelist
    \ifnum\currentp@ge=\f@nline
    \global\p@gelist=\expandafter\m@kelist\the\t@mplistX0
    \else\currentp@ge=\f@nline
    \global\p@gelist=\expandafter\m@kelist\the\t@mplistX1\fi\fi
    \let\next=\re@dfns\fi\next}
\def\symbols#1{\symb@ls={#1}\s@mb@ls=\symb@ls} 
\def\bigsymbol{\textstyle}
\symbols{\bigsymbol\ast,\dagger,\ddagger,\sharp,\flat,\natural,\star}
\def\ftnumbers{\ftn@mberstrue} \def\ftsymbols{\ftn@mbersfalse}
\def\paginal{\byp@ge} \def\resetftnumbers{\ftn@mber=1}
\def\ftnote#1{\defm@rk\expandafter\expandafter\expandafter\footnote
    \expandafter\b@@st\m@rk{#1}}

%Miscellaneous macros:
\long\def\jump#1\endjump{}
\def\ssum{\mathop{\lower .1em\hbox{$\textstyle\Sigma$}}\nolimits}

\def\qed{\nobreak\kern 1em \vrule height .5em width .5em depth 0em}
\def\newneq{\hbox{\rlap{\hbox to 1\wd9{\hss$=$\hss}}\raise .1em 
   \hbox to 1\wd9{\hss$\scriptscriptstyle/$\hss}}}
\def\subsetne{\setbox9 = \hbox{$\subset$}\mathrel{\hbox{\rlap
   {\lower .4em \newneq}\raise .13em \hbox{$\subset$}}}}
\def\supsetne{\setbox9 = \hbox{$\subset$}\mathrel{\hbox{\rlap
   {\lower .4em \newneq}\raise .13em \hbox{$\supset$}}}}

%Blackboard bold:
\def\vbar{\mathchoice{\vrule height6.3ptdepth-.5ptwidth.8pt\kern-.8pt}
   {\vrule height6.3ptdepth-.5ptwidth.8pt\kern-.8pt}
   {\vrule height4.1ptdepth-.35ptwidth.6pt\kern-.6pt}
   {\vrule height3.1ptdepth-.25ptwidth.5pt\kern-.5pt}}
\def\f@dge{\mathchoice{}{}{\mkern.5mu}{\mkern.8mu}}
\def\b@c#1#2{{\rm \mkern#2mu\vbar\mkern-#2mu#1}}
\def\b@b#1{{\rm I\mkern-3.5mu #1}}
\def\b@a#1#2{{\rm #1\mkern-#2mu\f@dge #1}}
\def\bb#1{{\count4=`#1 \advance\count4by-64 \ifcase\count4\or\b@a A{11.5}\or
   \b@b B\or\b@c C{5}\or\b@b D\or\b@b E\or\b@b F \or\b@c G{5}\or\b@b H\or
   \b@b I\or\b@c J{3}\or\b@b K\or\b@b L \or\b@b M\or\b@b N\or\b@c O{5} \or
   \b@b P\or\b@c Q{5}\or\b@b R\or\b@a S{8}\or\b@a T{10.5}\or\b@c U{5}\or
   \b@a V{12}\or\b@a W{16.5}\or\b@a X{11}\or\b@a Y{11.7}\or\b@a Z{7.5}\fi}}

\catcode`\X=11 \catcode`\@=12

%  *** end including mathdefs.tex *** 
% % \input citeadd
%  *** start including citeadd.tex *** 
%   citeadd -- a few additions for 
% files from alice that were procesed with "citealice"

\expandafter\ifx\csname citeadd.tex\endcsname\relax
\expandafter\gdef\csname citeadd.tex\endcsname{}
\else \message{Hey!  Apparently you were trying to
\string\input{citeadd.tex} twice.   This does not make sense.} 
\errmessage{Please edit your file (probably \jobname.tex) and remove
any duplicate ``\string\input'' lines}\endinput\fi

%  *** end including citeadd.tex *** 
\sectno=-1   % start with sect 0
\localtags
%\jjtags
\NoBlackBoxes
\define\mr{\medskip\roster}
\define\sn{\smallskip\noindent}
\define\mn{\medskip\noindent}
\define\bn{\bigskip\noindent}
\define\ub{\underbar}
\define\wilog{\text{without loss of generality}}
\define\ermn{\endroster\medskip\noindent}
\define\dbca{\dsize\bigcap}
\define\dbcu{\dsize\bigcup}
\define \nl{\newline}
\magnification=\magstep 1
\documentstyle{amsppt}
% % \input alice2000
%  *** start including alice2000.tex *** 
% This file should be inputted whenever we use amsppt.sty and 
% old tex.  
%  Here we redefine \subjclass (use 1991 instead of 2000, otherwise 
% the following definition comes directly from 
%% 
%%              `amsppt.sty', generated 
%% on <1997/2/2> with the docstrip utility (2.2i).
%% 
%% The original source files were:
%% 
%% amsppt.doc 
%%% ====================================================================
%%% @AMSTeX-style-file{
%%%   filename  = "amsppt.sty",
%%%   version   = "2.1h",
%%%   date      = "1997/02/02",
%%%   time      = "09:27:44 EST",
%%%   checksum  = "56844 3264 16617 137829",
%%%   author    = "American Mathematical Society",
%%%   address   = "PO Box 6248, Providence, RI 02940-6248, USA",
%%%   telephone = "401-455-4080 or (in the USA) 800-321-4AMS",

{    % the braces make the catcode-change local. 
\catcode`@11

\ifx\alicetwothousandloaded@\relax
  \endinput\else\global\let\alicetwothousandloaded@\relax\fi

\gdef\subjclass{\let\savedef@\subjclass
 \def\subjclass##1\endsubjclass{\let\subjclass\savedef@
   \toks@{\def\usualspace{{\rm\enspace}}\eightpoint}%
   \toks@@{##1\unskip.}%
   \edef\thesubjclass@{\the\toks@
     \frills@{{\noexpand\rm2000 {\noexpand\it Mathematics Subject
       Classification}.\noexpand\enspace}}%
     \the\toks@@}}%
  \nofrillscheck\subjclass}
} 

%  *** end including alice2000.tex *** 
\pageheight{8.5truein}
\topmatter
\title{On ultraproducts of Boolean Algebras and irr} \endtitle
\author {Saharon Shelah \thanks {\null\newline I would like to thank 
Alice Leonhardt for the beautiful typing. \null\newline
This research was supported by The Israel Science Foundation founded
by the Israel Academy of Sciences and Humanities. Publication 703.}
\endthanks} \endauthor  
\affil{Institute of Mathematics\\
 The Hebrew University\\
 Jerusalem, Israel
 \medskip
 Rutgers University\\
 Mathematics Department\\
 New Brunswick, NJ  USA} \endaffil
\endtopmatter
\document  
% % \input alice2jlem
%  *** start including alice2jlem.tex *** 
%% # Keywords  Input file to be used for texing Alice's files

\expandafter\ifx\csname alice2jlem.tex\endcsname\relax
  \expandafter\xdef\csname alice2jlem.tex\endcsname{\the\catcode`@}
\else \message{Hey!  Apparently you were trying to
\string\input{alice2jlem.tex}  twice.   This does not make sense.}
\errmessage{Please edit your file (probably \jobname.tex) and remove
any duplicate ``\string\input'' lines}\endinput\fi

% % \input bib4plain
%  *** start including bib4plain.tex *** 
\expandafter\ifx\csname bib4plain.tex\endcsname\relax
  \expandafter\gdef\csname bib4plain.tex\endcsname{}
\else \message{Hey!  Apparently you were trying to \string\input
  bib4plain.tex twice.   This does not make sense.}
\errmessage{Please edit your file (probably \jobname.tex) and remove
any duplicate ``\string\input'' lines}\endinput\fi

%  This file should be inputted if you want to use 
%  bibtex fom within plain TeX. 
      % Not really need for standard
       % bibtex files, but these commands
\def\renewcommand{\newcommand}	       % are used in our literal-unsrt.bst
\edef\cite{\the\catcode`@}%
\catcode`@ = 11
\let\@oldatcatcode = \cite
\chardef\@letter = 11
\chardef\@other = 12
%
%
% Next come some things that will be useful later.
%
% Make an outer definition into an inner one (due to Chris Thompson).
% The arguments should be the control sequence to be defined, and the
% new of the \outer control sequence, as characters; the control
% sequence #1 is defined to be just the same as \csname#2\endcsname, but
% not \outer.  For example, \@innerdef\innernewcount{newcount} would
% define \innernewcount to be a non-outer version of \newcount.
%
\def\@innerdef#1#2{\edef#1{\expandafter\noexpand\csname #2\endcsname}}%
%
% We use \@innerdef to make some of our allocations local, because
% Eplain includes our code inside a conditional.  We put @'s in the
% names to minimize the (already small) chance of conflicts.
%
\@innerdef\@innernewcount{newcount}%
\@innerdef\@innernewdimen{newdimen}%
\@innerdef\@innernewif{newif}%
\@innerdef\@innernewwrite{newwrite}%
%
%
% Swallow one parameter.
%
\def\@gobble#1{}%
%
%
% Use TeX 3.0's \inputlineno to get the line number, for better error
% messages, but if we're using an old version of TeX, don't do anything.
%
\ifx\inputlineno\@undefined
   \let\@linenumber = \empty % Pre-3.0.
\else
   \def\@linenumber{\the\inputlineno:\space}%
\fi
%
%
% The following macro \@futurenonspacelet (from the TeXbook) behaves
% essentially like \futurelet except that it discards any implicit or
% explicit space tokens that intervene before a nonspace is scanned:
%
\def\@futurenonspacelet#1{\def\cs{#1}%
   \afterassignment\@stepone\let\@nexttoken=
}%
\begingroup % The grouping here avoids stepping on an outside use of `\\'.
\def\\{\global\let\@stoken= }%
\\ % now \@stoken is a space token (\\ is a control symbol, so that
   % space after it is seen).
\endgroup
\def\@stepone{\expandafter\futurelet\cs\@steptwo}%
\def\@steptwo{\expandafter\ifx\cs\@stoken\let\@@next=\@stepthree
   \else\let\@@next=\@nexttoken\fi \@@next}%
\def\@stepthree{\afterassignment\@stepone\let\@@next= }%
%
%
% \@getoptionalarg\CS gets an optional argument from the input, enclosed
% in brackets, then expands \CS.  We set \@optionalarg to \empty if we
% don't find one, otherwise to the text of the argument.  This assumes
% the brackets don't have a funny category code.
%
\def\@getoptionalarg#1{%
   \let\@optionaltemp = #1%
   \let\@optionalnext = \relax
   \@futurenonspacelet\@optionalnext\@bracketcheck
}%
%
% The \expandafter's in this macro let us avoid the use of \aftergroup,
% which is somewhat more expensive.
%
\def\@bracketcheck{%
   \ifx [\@optionalnext
      \expandafter\@@getoptionalarg
   \else
      \let\@optionalarg = \empty
      % We can't do the \temp after the \fi, because then the \temp gets
      % in the way of reading the optional argument from the input, if
      % we do have one.
      \expandafter\@optionaltemp
   \fi
}%
\def\@@getoptionalarg[#1]{%
   \def\@optionalarg{#1}%
   \@optionaltemp
}%
%
%
% From LaTeX.
%
\def\@nnil{\@nil}%
\def\@fornoop#1\@@#2#3{}%
\def\@for#1:=#2\do#3{%
   \edef\@fortmp{#2}%
   \ifx\@fortmp\empty \else
      \expandafter\@forloop#2,\@nil,\@nil\@@#1{#3}%
   \fi
}%
\def\@forloop#1,#2,#3\@@#4#5{\def#4{#1}\ifx #4\@nnil \else
       #5\def#4{#2}\ifx #4\@nnil \else#5\@iforloop #3\@@#4{#5}\fi\fi
}%
\def\@iforloop#1,#2\@@#3#4{\def#3{#1}\ifx #3\@nnil
       \let\@nextwhile=\@fornoop \else
      #4\relax\let\@nextwhile=\@iforloop\fi\@nextwhile#2\@@#3{#4}%
}%
%
%
% This macro tests if a file \jobname.#1 exists, and sets \if@fileexists
% appropriately.  If an optional argument is given, it is used as the
% root part of the filename instead of \jobname.
%
\@innernewif\if@fileexists
\def\@testfileexistence{\@getoptionalarg\@finishtestfileexistence}%
\def\@finishtestfileexistence#1{%
   \begingroup
      \def\extension{#1}%
      \immediate\openin0 =
         \ifx\@optionalarg\empty\jobname\else\@optionalarg\fi
         \ifx\extension\empty \else .#1\fi
         \space
      \ifeof 0
         \global\@fileexistsfalse
      \else
         \global\@fileexiststrue
      \fi
      \immediate\closein0
   \endgroup
}%
%
%
%% [[[start of BibTeX-specific stuff]]]
%
% Now come the four main LaTeX commands and their associated .aux
% commands.  Just as in LaTeX, \bibliographystyle defines the BibTeX
% style name (.bst file, that is), and \bibliography defines the
% database (.bib) file(s).  The corresponding .aux-file commands are
% \bibstyle and \bibdata, which are there only for BibTeX's (but not
% LaTeX's) use.
%
\def\bibliographystyle#1{%
   \@readauxfile
   \@writeaux{\string\bibstyle{#1}}%
}%
\let\bibstyle = \@gobble
%
% As well as writing the \bibdata command to tell BibTeX which .bib
% files to read, we read the .bbl file that BibTeX (or a person,
% conceivably) has produced.  We use \bblfilebasename as the root of the
% filename to read; this defaults to \jobname.
%
\let\bblfilebasename = \jobname
\def\bibliography#1{%
   \@readauxfile
   \@writeaux{\string\bibdata{#1}}%
   \@testfileexistence[\bblfilebasename]{bbl}%
   \if@fileexists
      % We just output a non-discardable item (the `whatsit' with the
      % \bibdata command).  This means that the glue that will be
      % inserted next (\parskip or \baselineskip, most likely) will be a
      % legal breakpoint.  Most likely, this is after some kind of
      % heading, where we don't want to allow a page break.  So:
      \nobreak
      \@readbblfile
   \fi
}%
\let\bibdata = \@gobble
%
% The \nocite{label,label,...} command writes its argument to \@auxfile,
% unless instructed not to, but produces no text in the document.  Both
% the \nocite and \cite commands produce \citation commands in the .aux file.
%
\def\nocite#1{%
   \@readauxfile
   \@writeaux{\string\citation{#1}}%
}%
\@innernewif\if@notfirstcitation
%
% \cite[note]{label,label,...} produces the citations for the labels as
% well.  If the optional argument `note' is present, it's added after
% the labels.  Since \cite calls \nocite to do its .aux-file writing,
% \cite doesn't need to call \@readauxfile (\nocite does).
%
\def\cite{\@getoptionalarg\@cite}%
%
% Typeset the citations for the labels in #1, followed by the note, if
% it exists.  To change the citation's format in the text, redefine one
% or more `\print...' macros, whose defaults appear later in this file.
%
\def\@cite#1{%
   % Remember the optional argument, in case one of the macros we call
   % below ends up looking for an optional argument itself.  For
   % example, if a \cite[note] triggers reading the .aux file, then the
   % [note] would be clobbered, since \@testfileexistence looks for an
   % optional arg.
   \let\@citenotetext = \@optionalarg
   % Start printing the text, beginning with a left bracket by default.
   \printcitestart
   % It's complicated, but because \nocite puts a `whatsit' onto the list,
   % \nocite should follow \printcitestart.  It's conceivable, but very
   % unlikely, that this `whatsit' will cause a problem (glue that doesn't
   % disappear when you want it to is the most likely symptom), requiring
   % a change either to \printcitestart or to the label that the .bst file
   % produces.
   \nocite{#1}%
   \@notfirstcitationfalse
   \@for \@citation :=#1\do
   {%
      \expandafter\@onecitation\@citation\@@
   }%
   \ifx\empty\@citenotetext\else
      \printcitenote{\@citenotetext}%
   \fi
   \printcitefinish
}%
\def\@onecitation#1\@@{%
   \if@notfirstcitation
      \printbetweencitations
   \fi
   \expandafter \ifx \csname\@citelabel{#1}\endcsname \relax
      \if@citewarning
         \message{\@linenumber Undefined citation `#1'.}%
      \fi
      % Give it a dummy definition:
      \expandafter\gdef\csname\@citelabel{#1}\endcsname{%
% Change: marginal remark added, goldstrn@math.huji.ac.il, 
% goldstern@tuwien.ac.at, May 1996 mg
%  !!! change !!!
\strut
\vadjust{\vskip-\dp\strutbox
\vbox to 0pt{\vss\parindent0cm \leftskip=\hsize 
\advance\leftskip3mm
\advance\hsize 4cm\strut\openup-4pt 
\rightskip 0cm plus 1cm minus 0.5cm ?  #1 ?\strut}}
         {\tt
            \escapechar = -1
            \nobreak\hskip0pt
            \expandafter\string\csname#1\endcsname
            \nobreak\hskip0pt
         }%
      }%
   \fi
   % Now produce the text, whether it was undefined or not.
   \csname\@citelabel{#1}\endcsname
   \@notfirstcitationtrue
}%
%
% Given a label `foo', the macro `\b@foo' is supposed to
% hold the text that should be produced.
%
\def\@citelabel#1{b@#1}%
%
% So, how does a citation label get defined?  When we read the .bbl file
% (below), a \bibitem writes out a \@citedef command.  And when we read
% the \@citedef, we define \@citelabel{#1}, where #1 is the user's
% label.
%
\def\@citedef#1#2{\expandafter\gdef\csname\@citelabel{#1}\endcsname{#2}}%
%
%
% Reading the .bbl file also produces the typeset bibliography.  Please
% notice, however, that we do not produce the title for the references
% (e.g., `References'), as LaTeX does.  The formatting and spacing of
% that title, whether it should go into the headline, and so on, are all
% things determined by your format.  We cannot know those things in
% advance.  If you wish, you can define \bblhook to produce the title.
% Or just do it before the \bibliography command.
%
\def\@readbblfile{%
   % Define a counter to tell us which item number we are on, unless
   % we've already defined it (because the document has more than one
   % bibliography).
   \ifx\@itemnum\@undefined
      \@innernewcount\@itemnum
   \fi
   \begingroup
      \def\begin##1##2{%
         % ##1 is just `thebibliography'.
         % ##2 is the widest label.
         % We set (new dimen) \biblabelwidth based on the widest label
         \setbox0 = \hbox{\biblabelcontents{##2}}%
         \biblabelwidth = \wd0
      }%
      \def\end##1{}% ##1 is `thebibliography' again.
      %
      % Here we have two possibilities:
      % \bibitem[typesetlabel]{citationlabel}
      % \bibitem{citationlabel}
      % If we have the second of these, the citations are numbered, starting
      % from one; we use our own count register \@itemnum for this.
      %
      \@itemnum = 0
      \def\bibitem{\@getoptionalarg\@bibitem}%
      \def\@bibitem{%
         \ifx\@optionalarg\empty
            \expandafter\@numberedbibitem
         \else
            \expandafter\@alphabibitem
         \fi
      }%
      \def\@alphabibitem##1{%
         % Need \xdef here for various reasons.
         \expandafter \xdef\csname\@citelabel{##1}\endcsname {\@optionalarg}%
         % Left-justify alpha labels, unless \biblabel{pre,post}contents
         % are already defined.
         \ifx\biblabelprecontents\@undefined
            \let\biblabelprecontents = \relax
         \fi
         \ifx\biblabelpostcontents\@undefined
            \let\biblabelpostcontents = \hss
         \fi
         \@finishbibitem{##1}%
      }%
      \def\@numberedbibitem##1{%
         \advance\@itemnum by 1
         \expandafter \xdef\csname\@citelabel{##1}\endcsname{\number\@itemnum}%
         % Right-justify numeric labels, unless \biblabel{pre,post}contents
         % are already defined.
         \ifx\biblabelprecontents\@undefined
            \let\biblabelprecontents = \hss
         \fi
         \ifx\biblabelpostcontents\@undefined
            \let\biblabelpostcontents = \relax
         \fi
         \@finishbibitem{##1}%
      }%
      \def\@finishbibitem##1{%
         \biblabelprint{\csname\@citelabel{##1}\endcsname}%
         \@writeaux{\string\@citedef{##1}{\csname\@citelabel{##1}\endcsname}}%
         \ignorespaces
      }%
      %
      % Do the printing (we're producing the bibliography, remember).
      %
      \let\em = \bblem
      \let\newblock = \bblnewblock
      \let\sc = \bblsc
      % Punctuation won't affect spacing;
      \frenchspacing
      % the penalties below are from LaTeX's [article,book,report].sty;
      \clubpenalty = 4000 \widowpenalty = 4000
      % the next two values come from LaTeX's \sloppy command;
      \tolerance = 10000 \hfuzz = .5pt
      \everypar = {\hangindent = \biblabelwidth
                      \advance\hangindent by \biblabelextraspace}%
      \bblrm
      % the \parskip is a guess at what looks good;
      \parskip = 1.5ex plus .5ex minus .5ex
      % and the space between label and text comes from LaTeX's \labelsep.
      \biblabelextraspace = .5em
      \bblhook
      \input \bblfilebasename.bbl
   \endgroup
}%
%
% The widest label's width is useful for redefining \biblabelprint;
% you redefine \biblabelwidth, in effect, by redefining the
% \biblabelcontents macro that appears below.  And \biblabelextraspace,
% which is redefinable inside \bblhook, is added to \biblabelwidth to
% determine the amount of hanging indentation.
%
\@innernewdimen\biblabelwidth
\@innernewdimen\biblabelextraspace
%
% Now come the main macros that are related to the printing of the
% bibliography.  Since you might want to redefine them, they are given
% default definitions outside of \@readbblfile.
%
% The first one controls the printing of a bibliography entry's label.
% If you change it, make sure that it starts with something like
% \noindent or \indent or \leavevmode that puts TeX into horizontal mode
% (even if the label itself is empty); otherwise, the hanging
% indentation will get messed up in certain circumstances.
%
\def\biblabelprint#1{%
   \noindent
   \hbox to \biblabelwidth{%
      \biblabelprecontents
      \biblabelcontents{#1}%
      \biblabelpostcontents
   }%
   \kern\biblabelextraspace
}%
%
% If you are using numeric labels, and you want them left-justified
% (numeric labels by default are right-justified), do something like:
%     \def\biblabelprecontents{\relax}
%     \def\biblabelpostcontents{\hss}
%
% By default the labels are typeset in \bblrm, and enclosed in brackets.
%
\def\biblabelcontents#1{{\bblrm [#1]}}%
%
% The main text, too, is typeset using \bblrm, which is \rm by default.
%
\def\bblrm{\rm}%
%
% Emphasis for producing, e.g., titles, is done with \it by default.
%
\def\bblem{\it}%
%
% Some styles use a caps-and-small-caps font for author names.  LaTeX
% defines an \sc command but plain TeX doesn't, so we need one here.
% The definition below doesn't load the font unless it's needed, but it
% tries to load only the 10pt version, because it might not exist at
% other point sizes.
%
\def\bblsc{\ifx\@scfont\@undefined
              \font\@scfont = cmcsc10
           \fi
           \@scfont
}%
%
% The major parts of an entry are separated with \bblnewblock.  The
% numbers below are taken from LaTeX's `article' style.
%
\def\bblnewblock{\hskip .11em plus .33em minus .07em }%
%
% Here's where you stick any other bibliography-formatting goodies, or
% redefine the values above.
%
\let\bblhook = \empty
%
%
% Here are the four default definitions for formatting the in-text
% citations.  These are what you redefine (after your \input btxmac but
% before your \bibliography) to get parens instead of brackets, or
% superscripts, or footnotes, or whatever.
%
\def\printcitestart{[}%         left bracket
\def\printcitefinish{]}%        right bracket
\def\printbetweencitations{, }% comma, space
\def\printcitenote#1{, #1}%     comma, space, note (if it exists)
%
% That scheme is pretty flexible.  For example you could use
%     \def\printcitestart{\unskip $^\bgroup}
%     \def\printcitefinish{\egroup$}
%     \def\printbetweencitations{,}
%     \def\printcitenote#1{\hbox{\sevenrm\space (#1)}}
%     \font\eighttt = cmtt8
%     \scriptfont\ttfam = \eighttt
% to get superscripted in-text citations.  (The scriptfont stuff
% exists only to print an undefined citation; it's in cmtt8 because
% there is no cmtt7.)  To get something radically different, however,
% you'll have to define your own \cite command.
%
% When we read `\citation' from the .aux file, it means nothing.
%
\let\citation = \@gobble
%
%
% Now comes the stuff for dealing with LaTeX's \newcommand.  As
% mentioned earlier, this \newcommand will redefine a preexisting
% command; that's different from how LaTeX's \newcommand behaves.
%
\@innernewcount\@numparams
%
% \newcommand{\foo}[n]{text} defines the control sequence \foo to have
% n parameters, and replacement text `text'.
%
\def\newcommand#1{%
   \def\@commandname{#1}%
   \@getoptionalarg\@continuenewcommand
}%
%
% Figure out if this definition has parameters.
%
\def\@continuenewcommand{%
   % If no optional argument, we have zero parameters.  Otherwise, we
   % have that many.
   \@numparams = \ifx\@optionalarg\empty 0\else\@optionalarg \fi \relax
   \@newcommand
}%
%
% \@numparams is how many arguments this command has.  The name of the
% command is \@commandname.  The replacement text for the new macro is #1.
%
\def\@newcommand#1{%
   \def\@startdef{\expandafter\edef\@commandname}%
   \ifnum\@numparams=0
      \let\@paramdef = \empty
   \else
      \ifnum\@numparams>9
         \errmessage{\the\@numparams\space is too many parameters}%
      \else
         \ifnum\@numparams<0
            \errmessage{\the\@numparams\space is too few parameters}%
         \else
            \edef\@paramdef{%
               % This is disgusting, but \loop doesn't work inside \edef,
               % because \body isn't defined.
               \ifcase\@numparams
                  \empty  No arguments.
               \or ####1%
               \or ####1####2%
               \or ####1####2####3%
               \or ####1####2####3####4%
               \or ####1####2####3####4####5%
               \or ####1####2####3####4####5####6%
               \or ####1####2####3####4####5####6####7%
               \or ####1####2####3####4####5####6####7####8%
               \or ####1####2####3####4####5####6####7####8####9%
               \fi
            }%
         \fi
      \fi
   \fi
   \expandafter\@startdef\@paramdef{#1}%
}%
%
%% [[[end of BibTeX-specific stuff]]]
%
%
% Names of references (arguments given in the \cite and \nocite
% commands) and file names (arguments given in the \bibliography and
% \bibliographystyle commands) are recorded in \jobname.aux, called the
% \@auxfile in these macros.  Here's how they get read in.
%
\def\@readauxfile{%
   \if@auxfiledone \else % remember: \@auxfiledonetrue if \noauxfile is defined
      \global\@auxfiledonetrue
      \@testfileexistence{aux}%
      \if@fileexists
         \begingroup
            % Because we might be in horizontal mode when \@readauxfile
            % is called (if it's in response to a \cite or \nocite), we
            % want to ignore all the would-be spaces at the ends of
            % lines in the aux file.  Fortunately, it's highly unlikely
            % an end-of-line might actually be desired.
            % And because we don't change the category code of anything
            % but @, primitives like \gdef can't be used to define labels
            % in the aux file.  The solution adopted by btxmac.tex is to
            % write `\@citedef{LABEL}{DEFINITION}' to the aux file, and
            % use \csname on LABEL.
            \endlinechar = -1
            \catcode`@ = 11
            \input \jobname.aux
         \endgroup
      \else
         \message{\@undefinedmessage}%
         \global\@citewarningfalse
      \fi
      \immediate\openout\@auxfile = \jobname.aux
   \fi
}%
%
% The \@readauxfile macro does all that work the first time it's called.
% Since it's called once for every \cite, \nocite, \bibliography, and
% \bibliographystyle command that the user issues, we need to remember
% whether the work's been done.  It's considered done if we're not to do
% it---that is, if \noauxfile is defined.
%
\newif\if@auxfiledone
\ifx\noauxfile\@undefined \else \@auxfiledonetrue\fi
%
% It's conceivable you'd want to change how other characters are read;
% to do that, change their category code before doing \input btxmac.
%
%
% After reading the .aux file, \@readauxfile opens it for writing.
% The \@writeaux macro does the actual writing (as long as
% \noauxfile is undefined).
%
\@innernewwrite\@auxfile
\def\@writeaux#1{\ifx\noauxfile\@undefined \write\@auxfile{#1}\fi}%
%
%
% A macro package that uses btxmac.tex might define
% \@undefinedmessage (before doing an \input btxmac).
%
\ifx\@undefinedmessage\@undefined
   \def\@undefinedmessage{No .aux file; I won't give you warnings about
                          undefined citations.}%
\fi
%
% Even if citations are undefined, we want to complain only if
% \@citewarningtrue.  The default is to set \@citewarningtrue unless
% \noauxfile is defined.  Again, a macro package that uses
% btxmac.tex might want to redefine this.
%
\@innernewif\if@citewarning
\ifx\noauxfile\@undefined \@citewarningtrue\fi
%
%
% Finally, before leaving we restore @'s old category code.
%
\catcode`@ = \@oldatcatcode

%  *** end including bib4plain.tex *** 
  % This will define \cite and make sure it works as in latex

\def\widestnumber#1#2{}
  % Our amstex-ppt style does not know about \widestnumber

\def\rm{\fam0 \tenrm}

\def\fakesubhead#1\endsubhead{\bigskip\noindent{\bf#1}\par}

% % \input rsfs
%  *** start including rsfs.tex *** 

% # Keywords: Script or Calligraphic (Caligraphic) letters with the RSFS Font

% The story so far:    July 1998 -- Saharon would like to have a
% ``nicer'' calligraphic font. In particualr, the leters S and P in
% the usual calligraphic font do not look ``special'' enough. 
% 
% I found out that ``rsfs'' (``Ralph Smith Formal Script'') may be
% what he wants.   I installed the mf file, the .tfm file, as well as
% a few pk files in ~/TeX/rsfs.    Let's hope that this is enough.
% Using amstex, all you have to do is to \input rsfs.tex 
% Files prepared with citealice willdothis automatically. 
%
%  Note:  for some reason xdvi calls MakeTeXpk, then Maketexpk
%  complains about wrong resolution, but still writes commands to
%  missfont.log...  
%

% we redefine a macro inside amstex's \Cal command , so that it calls
% our nice font ``rsfs'' rather than the usual calligraphic font. 
% Note thisworks for amstex only.   
% In plain tex, would have to add definitions of \Cal
% in latex... we should insteaduse mathrsfs.sty
% 

\font\textrsfs=rsfs10
\font\scriptrsfs=rsfs7
\font\scriptscriptrsfs=rsfs5

\newfam\rsfsfam
\textfont\rsfsfam=\textrsfs
\scriptfont\rsfsfam=\scriptrsfs
\scriptscriptfont\rsfsfam=\scriptscriptrsfs

\edef\oldcatcodeofat{\the\catcode`\@}
\catcode`\@11

\def\Cal@@#1{\noaccents@ \fam \rsfsfam #1}

\catcode`\@\oldcatcodeofat

%  *** end including rsfs.tex *** 

\expandafter\ifx \csname margininit\endcsname \relax\else\margininit\fi

%  *** end including alice2jlem.tex *** 
\bigskip

\head {Anotated Content} \endhead  \resetall 
 % \resetall 
\bn
\S1 $\quad$ Consistent inequality
\mr
\item "{{}}"  [We prove the consistency of 
irr$(\dsize \prod_{i < \kappa} B_i/D) < \dsize \prod_{i < \kappa}$ 
irr$(B_i)/D$ where $D$ is an ultrafilter on $\kappa$ and each $B_i$ is a
Boolean Algebra.  This solves the last problem of this form from the Monk's
list of problems in \cite{M2}, that is number 35.  
The solution applies to many other properties, e.g. Souslinity.]
\endroster
\bn
\S2 $\quad$ Consistency for small cardinals
\mr
\item "{{}}"  [We get similar results with $\kappa = \aleph_1$ (easily
we cannot have it for $\kappa = \aleph_0$) and Boolean Algebras $B_i$
($i < \kappa$) of cardinality $< \beth_{\omega_1}$.]
\ermn
This article continues Magidor Shelah \cite{MgSh:433} and Shelah Spinas
\cite{ShSi:677}, but does not rely on them: see \cite{M2} on background.
\newpage

\head {\S1 Consistent inequality} \endhead  \resetall \sectno=1
 % \resetall 
\bigskip

\definition{\stag{1.1} Definition}  Assume $\mu < \lambda,\lambda$ is strongly
inaccessible Mahlo.  Let $B^* = B_\lambda$ be the Boolean Algebra freely
generated by $\{x_\alpha:\alpha < \lambda\}$ and for $u \subseteq \lambda$
let $B_u$ be the subalgebra of $B^*$ generated by $\{x_\alpha:\alpha \in u\}$.
\nl
1) We define a forcing notion $\Bbb Q = \Bbb Q^1_{\mu,\lambda}$ as follows:

$p \in \Bbb Q$ iff: for some $w^p = w[p]$ we have:
\mr
\widestnumber\item{$(iii)$}
\item "{$(i)$}"  $w^p  = w[p] \subseteq \lambda$
\sn
\item "{$(ii)$}"  $B^p = B[p]$ is a Boolean Algebra of the form
$B_{w[p]}/I^p$ where $I^p = I[p]$ is an ideal of $B_{w[p]}$ so $B^p$ is
generated by $\{x_\alpha/I:\alpha \in w^p\}$
\sn
\item "{$(iii)$}"  $x_\alpha/I \notin \langle \{x_\beta/I:\beta \in w^p \cap
\alpha\} \rangle_{B[p]}$, equivalently $x_\alpha \notin \langle \{x_\beta:\beta \in
w^p \cap \alpha\} \cup I \rangle_{B_{w[p]}}$
\sn
\item "{$(iv)$}"  for every strongly inaccessible $\chi \in (\mu,\lambda]$
we have $|w^p \cap \chi| < \chi$.
\ermn
The order is $p \le q$ iff $w^p \subseteq w^q$ and $I^p = I^q \cap B_{w[q]}$,
so abusing notation we think $B^p \subseteq B^q$, not distinguishing sometime
$x_\alpha$ from $x_\alpha/I \in B^p$ or (see below) from 
$x_\alpha/\underset\tilde {}\to I$ in
$\underset\tilde {}\to B$. \nl
2) We define $\underset\tilde {}\to I = \cup\{I^p:p \in
{\underset\tilde {}\to G_{{\Bbb Q}^1_{\mu,\kappa}}}\}$ and
$\underset\tilde {}\to B$ is defined as $B_\lambda/\underset\tilde {}\to I$.
\enddefinition
\bigskip

\proclaim{\stag{1.2} Claim}  For $\mu < \lambda$ as in Definition
\scite{1.1}, the forcing notion $\Bbb Q^1_{\mu,\lambda}$ is $\mu^+$-complete
(hence, add no new subsets to $\mu$), has cardinality $\lambda$, satisfies 
the $\lambda$-c.c., collapse no cardinal, change no cofinality, so
cardinal arithmetic is clear.
\endproclaim
\bigskip

\demo{Proof}  Like Easton forcing.
\enddemo
\bigskip

\proclaim{\stag{1.3} Claim}  For the forcing $\Bbb Q = \Bbb
Q^1_{\mu,\lambda}$ with $\mu,\lambda$ as in Definition \scite{1.1} we have \nl
1) $\Vdash_{\Bbb Q} ``\underset\tilde {}\to B$ is a
Boolean Algebra generated by $\{x_\alpha:\alpha < \lambda\}$ such that
$\alpha < \lambda \Rightarrow x_\alpha \notin \langle \{x_\beta:\beta <
\alpha\} \rangle_{\underset\tilde {}\to B}$, so $|\underset\tilde {}\to B| =
\lambda"$. \nl
2) $\Vdash_{\Bbb Q} ``\text{irr}^+(\underset\tilde {}\to B) = \lambda = irr
(\underset\tilde {}\to B)"$, see definition \scite{1.4} below. \nl
3) $\Vdash_{\Bbb Q} ``\text{if } y_\beta \in \underset\tilde {}\to B$ for $\beta <
\lambda$ \ub{then} for some $\beta_0 < \beta_1 < \beta_2 < \lambda$ we have
$\underset\tilde {}\to B \models y_{\beta_1} \cap y_{\beta_2} =
y_{\beta_0}"$. \nl
4) Let $B^*$ be a finite Boolean Algebra generated by $\{a^*,b^*,y^*_0,
\dotsc,y^*_{n(*)}\}$ such that $y^*_m \notin \langle \{y^*_\ell:\ell < m\}
\cup \{a^*,b^*\} \rangle,0 < a^* < y^*_\ell < b^* < 1$.
\mr
\item "{${{}}$}"  \ub{Then} it is forced, 
$(\Vdash_{{\Bbb Q}^1_{\mu,\lambda}})$ that: \nl
if $y_\beta \in \underset\tilde {}\to B$ for $\beta < \lambda$ and 
$\beta \ne \gamma \Rightarrow y_\beta \ne y_\gamma$
\ub{then} we can find $a,b$ in $\underset\tilde {}\to B$ satisfying
$0 < a < b < 1$ and
$\beta_0 < \ldots < \beta_{n(*)} < \lambda$ such that
{\roster
\itemitem{$ (\alpha)$ }   $\underset\tilde {}\to B \models ``a < y_{\beta_\ell} < b"$
\sn
\itemitem{ $(\beta)$ }  there is an embedding $f$ of $B^*$ into
$\underset\tilde {}\to B$ mapping $a^*$ to $a,b^*$ to $b$ and $y_\ell$
to $y^*_{\beta_\ell}$ for $\ell = 0,\dotsc,n(*)$.
\endroster}
\endroster
\endproclaim
\bn
Recalling
\definition{\stag{1.4} Definition}  For a Boolean Algebra $B$ let: \nl
1) $X \subseteq B$ is called irredundant, if no $x \in X$ belongs to the
subalgebra $\langle X \backslash \{x\} \rangle_B$ of $B$ generated by
$X \backslash \{x\}$. \nl
2) irr$^+(B) = \cup\{|X|^+:X \subseteq B \text{ is irredundent}\}$. \nl
3) irr$(B) = \cup\{|X|:X \subseteq B \text{ is irredundent}\}$ so
irr$(B)$ is irr$^+(B)$ if the latter is a limit cardinal and is the
predecessor of irr$^+(B)$ if the later is a successor cardinal.
\enddefinition
\bigskip

\remark{Remark}  Concering \scite{1.3}  
on the case $\kappa = \aleph_1$, see Rubin \cite{Ru83}, generally see
\cite{Sh:128}, \cite{Sh:e}.
\endremark
\bigskip

\demo{Proof of 1.3}  1) Should be clear. \nl
2) Clearly for every $\chi < \lambda$ and $p \in \Bbb Q^1_{\mu,\lambda}$ we 
can find $\alpha < \lambda$ such that $\alpha > \chi$ and 
$w^p \cap [\alpha,\alpha + \chi) =
\emptyset$ hence we can find $q$ such that $p \le q \in 
\Bbb Q^1_{\mu,\lambda}$ and $w^q = w^p \cup [\alpha,\alpha + \chi)$ and in
$B^q$ the set $\{x_\beta:\beta \in [\alpha,\alpha + \chi)\}$ is independent, 
hence
$q \Vdash ``\text{irr}^+(\underset\tilde {}\to B) > \chi"$.  So we get
$\Vdash ``\text{irr}^+(B) \ge \lambda$.  To prove equality use part (3). \nl
3) Assume toward contradiction that $p \Vdash ``\langle
{\underset\tilde {}\to y_\beta}:\beta < \lambda \rangle$ is a counterexample".
We can find for each $\beta < \lambda$ a quadruple 
$(p_\beta,n_\beta,\langle \alpha_{\beta,\ell}:
\ell < n_\beta \rangle,\sigma_\beta)$ such that:
\mr
\widestnumber\item{$(iii)$}
\item "{$(i)$}"  $p \le p_\beta \in \Bbb Q^1_{\mu,\lambda}$
\sn
\item "{$(ii)$}"  $n_\beta < \omega$
\sn
\item "{$(iii)$}"  $\alpha_{\beta,\ell} \in w^{p_\beta}$ increasing with
$\ell$
\sn
\item "{$(iv)$}"  $\sigma_\beta(x_0,\dotsc,x_{n_\beta-1})$ is
a Boolean term
\sn
\item "{$(v)$}"  $p_\beta \Vdash ``\text{in } \underset\tilde {}\to B \text{ we have }
{\underset\tilde {}\to y_\beta} = \sigma_\beta(x_{\alpha_{\beta,0}},
x_{\alpha_{\beta,1}},\dotsc, x_{\alpha_{\beta,n_\beta-1}})"$ call the 
latter $y_\beta$, so by part (1) \wilog \, $\{\alpha_{\beta,\ell}:\ell
< n_\beta\} \subseteq w^{p_\beta}$ hence $y_\beta$ is a member of $B_{w[p]}$.
\ermn
So we can choose a stationary $S \subseteq \{\chi:\chi \text{ strongly
inaccessible}, \mu < \chi < \lambda\}$ and 
$n,\sigma,m,\langle \alpha_\ell:\ell < m \rangle,w,r$ \ub{such that} 
for every $\beta \in S$ we have:
$n_\beta = n \and \sigma_\beta = \sigma,\ell < m \Rightarrow 
\alpha_{\beta,\ell} = \alpha_\ell,\ell \in [m,n) 
\Rightarrow \alpha_{\beta,\ell} \ge \beta$ and $w^{p_\beta} \cap
\beta = w$.
Without loss of generality also $\alpha < \beta
\in S \Rightarrow w^{p_\alpha} \subseteq \beta$.  Without loss of generality
for $\beta_0,\beta_1$ in $S$ the mapping $F_{\beta_0,\beta_1} =
\text{ id}_w \cup \{\langle
(\alpha_{\beta_0,\ell},\alpha_{\beta_1,\ell}):\ell < n \rangle\}$ induces an
isomorphism $g_{\beta_1,\beta_0}$ from the Boolean Algebra 
$\langle \{x_\gamma:\gamma \in w\} \cup
\{\alpha_{\beta_{0,\ell}}:\ell < n\} \rangle_{B[p_{\beta_0]}}$ 
onto the Boolean
Algebra $\langle\{x_\gamma:\gamma \in w\} \cup \{x_{\beta_1,\ell}:\ell < n\} 
\rangle_{B[p_{\beta_1}]}$ that is $g_{\beta_1,\beta_0}$ maps $x_\gamma$ to
$x_\gamma$ for $\gamma \in w$ and maps $x_{\beta_0,\ell}$ to $x_{\beta_1,\ell}$
for $\ell < n$.  Choose $\beta_0 < \beta_1 < \beta_0$ and we define
$q \in Q^1_{\mu,\lambda}$ such that $w^q = w[p_{\beta_0}] \cup w[p_{\beta_1}]
\cup w[p_{\beta_2}]$ and $B^q$ is the
Boolean Algebra generated by $\{x_\alpha:\alpha \in w[p_{\beta_0}] \cup
w[p_{\beta_1}] \cup w[p_{\beta_2}]\}$ freely except the equations which 
hold in $p_{\beta_\ell}$ 
for each $\ell = 0,1,2$ and the equation $y_{\beta_1} \cap y_{\beta_2} =
y_{\beta_0}$, in other words $I^q$ is the ideal of $B_{w^q}$ generated
by $I[p_{\beta_0}] \cup I[p{\beta_1}] \cup
I[p_{\beta_2}] \cup \{y_{\beta_1} \cap y_{\beta_2} - y_{\beta_0},y_{\beta_0}
- y_{\beta_1} \cap y_{\beta_2}\}$.  
We should prove that $q \in Q^1_{\mu,\lambda}$ and 
$I[q] \cap B[p_{\beta_\ell}] = I[p_{\beta_\ell}]$ for $\ell = 0,1,2$
(the rest: $p_{\beta_\ell} \le q$ hence $p \le q$ and $q \Vdash ``
{\underset\tilde {}\to y_{\beta_\ell}} = y_{\beta_\ell}$ for $\ell = 0,1,2$ 
and
$y_{\beta_1} \cap y_{\beta_2} = y_{\beta_0}"$ should be clear).  Let $B_0$
be the trivial Boolean Algebra $\{0,1\}$.
\nl
For $w \subseteq \lambda$ and $f \in {}^w 2$ let $\hat f$ be the unique
homomorphism from the Boolean Algebra $B_w$ 
freely generated by $\{x_\alpha:\alpha
\in w\}$ to $\{0,1\}$ such that $\alpha \in w \Rightarrow \hat f(x_\alpha)
= f(\alpha)$.  For $p^* \in Q^1_{\mu,\lambda}$ let
${\Cal F}[p^*] = \{f:f \in {}^{(w^{p^*})}2$ and $\{x_\alpha:f(\alpha) = 1\}
\cup \{-x_\alpha:f(\alpha) =0\}$ generates an ultrafilter of $B[p^*]$.  
For each $f \in {\Cal F}[p^*]$ let 
$f^{[p^*]}$ be the homomorphism from $B[p^*]$ to $B_0$ induced
by $f$, i.e. $f^{[p^*]}(x_\alpha) = f(\alpha)$ for every $\alpha \in w$.
Clearly ${\Cal F}[p^*]$ gives all the information on $p^*$.  Define
$u = w^{p_{\beta_0}} \bigcup w^{p_{\beta_1}} \bigcup w^{p_{\beta_2}}$ and

$$
\align
{\Cal F} = \bigl\{ f:&f \in {}^u 2, \text{ and } \ell \le 2 \Rightarrow f
\restriction w[p_{\beta_\ell}] \in {\Cal F}[p_{\beta_\ell}] \text{ and}\\
  &B_0 \models ``\hat f(\sigma(\langle x_{\beta_1,\ell}:\ell < n \rangle)) 
\cap \hat f(\sigma(\langle x_{\beta_2,\ell}:\ell < n \rangle)) \\ 
  &= \hat f(\sigma(\langle x_{\beta_0,\ell}:\ell < n \rangle))" \bigr\}.
\endalign
$$
\mn
We need to show that ${\Cal F}$ is rich enough, clearly $\otimes_1 + \otimes_2 +
\otimes_2$ below suffice
\mr
\item "{$\bigotimes_1$}"  if $\ell \in \{0,1,2\}$ and $f_\ell \in {\Cal F}
[p_{\beta_\ell}]$ then there is $f \in {\Cal F}$ extending $f_\ell$.
\ermn
[Why?  For $m=0,1,2$ let 
$p'_{\beta_m}$ be the subalgebra of $B[p_{\beta_m}]$ generated by
$\{x_\gamma:\gamma \in w[p_{\beta_m}]$ and $\gamma < \beta_m \vee 
\gamma \in \{\alpha_{\beta_m,0},
\dotsc,\alpha_{\beta_m,n-1}\}\}$.  We define a homomorphism $h_\ell$ from
$p'_{\beta_\ell}$ to $B_0$ as $f^{[p_{\beta_\ell}]}_\ell \restriction B[p'_{\beta_\ell}]$ 
and define for $m=0,1,2$ a homomorphism $g_m$ from $B[p'_{\beta_m}]$ to 
$B_0$ such that: $\gamma \in w \Rightarrow g_m(x_\gamma) =
f_\ell(\gamma)$ and $\gamma = \beta_{m,k} \Rightarrow g_m(x_\gamma)
= f_\ell(\beta_{\ell,k})$.   The definitions are
compatible and let $h_m$ be $h_\ell$ if $\ell = m$ and any 
homomorphism from $B[p_{\beta_m}]$ to 
$B_0$ extending $g_m$ if $m \in \{0,1,2\} \backslash
\{\ell\}$, clearly exist.  Let $f_m \in {}^{w[p_{\beta_\ell}]} 2$ for
$m = 0,1,2$ be $f_m(\gamma) = h_m(x_\gamma)$; for $m=\ell$ the
definitions are compatible.  Lastly let 
$f = f_0 \cup f_1 \cup f_2$, easily $f_\ell \subseteq
f \in {\Cal F}$.]
\mr
\item "{$\bigotimes_2$}"  if $\ell \in \{0,1,2\},\alpha \in w[p_{\beta_\ell}]$
then there are $f',f'' \in {\Cal F}$ such that $f'(\alpha) \ne f''
(\alpha)$ but $f' \restriction (\alpha \cap u) = f'' \restriction (\alpha
\cap u)$.
\ermn
[Why?  As $p_{\beta_\ell} \in \Bbb Q^1_{\mu,\lambda}$ we can find $f'_\ell,
f''_\ell \in {\Cal F}[p_{\beta_\ell}]$ such that $f'_\ell(\alpha) \ne
f''_\ell(\alpha)$ but $f'_\ell \restriction (\alpha \cap w[p_{\beta_\ell}])
= f'' \restriction (\alpha \cap w[p_{\beta_\ell}])$.  Now
for $m \in \{0,1,2,\} \backslash \{\ell\}$ let $f'_m \in {\Cal F}
[p_{\beta_m}]$ extends $f'_\ell \circ F_{\beta_\ell,\beta_m}$ and 
$f''_m \in {\Cal F}[p_{\beta_m}]$ extends 
$f''_m \circ F_{\beta_\ell,\beta_m}$  both times as in the proof of $\otimes_1$.  
If $\ell = 0$, let $f' = f'_0 \cup f'_1 \cup f'_2 \in {\Cal F}$ and let
$f''= f''_0 \cup f''_1 \cup f''_2 \in {\Cal F}$ and 
we are done.  
Also if $\alpha < \beta_\ell$
(so $\alpha \in \dbca_{m \le 2} w[p_{\beta_m}]$) the same proof works.
So assume $\ell \ne 0,\alpha \notin \dbca_{m \le 2} w[p_{\beta_m}]$.  If
$(f'_\ell)^{[p_{\beta_\ell}]}(y_{\beta_\ell}) =
(f''_\ell)^{[p_{\beta_\ell}]}(y_{\beta_\ell})$ let $f' = f'_0 \cup f'_1
\cup f'_2,f'' = f''_\ell \cup (f' \upharpoonright (w[p_{\beta_0}] \cup
w[p_{\beta_{3 - \ell}}]))$, clearly O.K.
So \wilog \, assume
$(f'_\ell)^{[p_{\beta_\ell}]}(y_{\beta_\ell}) = 0,
(f''_\ell)^{[p_{\beta_\ell}]}(y_{\beta_\ell}) =1,\ell \in \{1,2\}$ and
$\alpha \in w[p_{\beta_\ell}] \backslash w[p_{\beta_0}]$; and then choose
$f' = f'_0 \cup f'_1 \cup f'_2$ as above and 
$f'' = f''_\ell \cup (f' \restriction (w[p_{\beta_0}]
\cup w[\beta_{\beta_{3 - \ell}}]))$.  Now check; the main point is that as
$\hat f'_{3 - \ell}(y_{\beta_{3 - \ell}}) = \hat f'_0(y_{\beta_0})$ 
we have $B_0 \models ``\hat f''(y_{\beta_1}) \cap \hat f''(y_{\beta_2}) = 
\hat f''(y_{\beta_\ell}) \cap \hat f''(y_{\beta_{3 -\ell}}) = 
\hat f''_\ell(y_{\beta_\ell}) \cap \hat f'(y_{\beta_{3 - \ell}})
= 1_{B_0} \cap \hat f'_{3 - \ell}(y_{\beta_{3 - \ell}}) = \hat f'_{3 - \ell}
(y_{\beta_{3 - \ell}}) = \hat f'_0(y_{\beta_0}) = \hat f''(y_{\beta_0})$". 
\nl
4) Similar proof (with $a,b$ now in $p_{\beta_\ell} \restriction \beta_\ell!$).
\hfill$\square_{\scite{1.3}}$\margincite{1.3}
\enddemo
\bigskip

\proclaim{\stag{1.5} Claim}  1) If $\Bbb Q = \Bbb Q^1_{\mu,\lambda} *
{\underset\tilde {}\to {\Bbb Q}^2}$ and $\Vdash_{{\Bbb Q}^1_{\mu,\lambda}} ``
{\underset\tilde {}\to {\Bbb Q}^2}$ satisfies the $(\lambda,3)$-Knaster 
condition (see below)", \ub{then} $\Vdash_{\underset\tilde {}\to {\Bbb Q}} ``
\text{irr}^+(\underset\tilde {}\to B) = \lambda"$. \nl
2) For \scite{1.3}(4), ``$(\lambda,n^* +1)$-Knaster" suffice to preserve the
condition. \nl
3) In part (1) we even get the conclusion of Claim \scite{1.3}(3).
\endproclaim
\bigskip

\definition{\stag{1.6} Definition}  1) The $\lambda$-Knaster condition 
says that among any $\lambda$ members there is a set of $\lambda$ which 
are pairwise compatible.  Recall that it is preserved by composition. \nl
2) For $n^* \le \omega$, the $(\lambda,n^* +1)$-Knaster condition says that among
any $\lambda$ member there is a set of $\lambda$ such that any $< 1 + n^*$
of them has a common upper bound.
\enddefinition
\bigskip

\demo{Proof of \scite{1.5}}  1), 3)  Clearly it suffices to prove
(3). \nl
Straight by \scite{1.4}(3), in fact, just such $\Bbb Q^2$
preserves the properties mentioned there in \scite{1.5}.  \nl
2) Similarly using \scite{1.4}(4).   \hfill$\square_{\scite{1.5}}$\margincite{1.5}
\enddemo
\bigskip

\proclaim{\stag{1.7} Theorem}  Suppose
\mr
\item "{$(a)$}"  $\bold V$ satisfies GCH above $\mu$ (for simplicity)
\sn
\item "{$(b)$}"  $\kappa$ is measurable, $\kappa < \chi < \mu$
\sn
\item "{$(c)$}"  $\mu$ is supercompact, Laver indestructible, in fact,
{\roster
\itemitem{ $(*)$ }   for
some $h_\ell:\mu \rightarrow {\Cal H}(\mu)$, (for $\ell = 0,1$) we have \ub{for every}
$(< \mu)$-directed complete forcing $\Bbb Q$, cardinal $\theta \ge \mu$ and $\Bbb Q$-name
$\underset\tilde {}\to x$ of a subset of $\theta$, \ub{there is} in $\bold V[G_{\Bbb Q}]$
a normal ultrafilter ${\Cal D}$ on $[\theta]^{< \mu}$ such that \nl
$\dsize \prod_{a \in [\theta]^{< \mu}}(h_1(a \cap \mu),h_2(a \cap \mu)) \cong 
(\theta,\underset\tilde {}\to x[G_{\Bbb Q}])$ 
\endroster}
\sn
\item "{$(d)$}"  $\lambda > \mu$ is strongly inaccessible, Mahlo and
$\lambda^*$ such that $\lambda^* = (\lambda^*)^\mu \ge \lambda$
\sn
\item "{$(e)$}"  $D^*$ is a normal ultrafilter on $\kappa$.
\ermn
\ub{Then} for some forcing notion $\Bbb P$ we have, in $\bold V^{\Bbb P}$:
\mr
\item "{$(\alpha)$}"  forcing with $\Bbb P$ collapse no cardinal of
$\bold V$ except those in the interval $(\mu^+,\lambda)$
\sn
\item "{$(\beta)$}"  forcing with $\Bbb P$ add no subsets to $\chi$
\sn
\item "{$(\gamma)$}"  $\mu$ is strong limit of cofinality $\kappa$ and 
$\langle \mu_i:i < \kappa \rangle$ is an increasing continuous
sequence of strong limit cardinals with limit $\mu$
\sn
\item "{$(\delta)$}"  for each $i < \kappa,\mu_i < \lambda_i \le \lambda^*_i
= (\lambda^*_i)^{\mu_i} = 2^{\mu_i}$ and we let
$\mu_\kappa = \mu,\lambda_\kappa = \lambda,\lambda^*_\kappa = \lambda^*$
\sn
\item "{$(\varepsilon)$}"  for each $i \le \kappa$ we have: $B_i$ is a Boolean
Algebra of cardinality $\lambda_i$ and irr$^+(B_i) = \lambda_i$
\sn
\item "{$(\zeta)$}"  for $i < \kappa,\lambda_i$ is a Mahlo cardinal
even strongly inaccessible,  but 
\sn
\item "{$(\eta)$}"  $\lambda = \lambda_\kappa$ is $\mu^{++}$
(this in $V^{\Bbb P}$)
\sn
\item "{$(\theta)$}"  $B = B_\kappa$ is 
isomorphic to $\dsize \prod_{i < \kappa} B_i/D^*$, hence
{\roster
\itemitem{ $\boxtimes$ }  {\rm irr}$^+(B) = \lambda = \mu^{++}$ so
{\rm irr}$(B) = \mu^+$ whereas {\rm irr}$(B_i) = \text{ {\rm irr\/}}^+(B_i)
= \lambda_i$ and $\dsize \prod_{i < \kappa} \lambda_i/D^* = \lambda$, so
{\rm irr}$(\dsize \prod_{i < \kappa} B_i/D^*) < \dsize \prod_{i < \kappa}
\text{ {\rm irr\/}}(B_i)/D^*$.
\endroster}
\endroster
\endproclaim
\bigskip

\demo{Proof}  Let $\Bbb Q_1 = \Bbb Q^1_{\mu,\lambda}$ and
$\underset\tilde {}\to B$ be from \scite{1.2}, let
$\Bbb Q_2$ be $\{f:f$ a partial function from $\lambda^*$ to $\{0,1\}$ with
domain of cardinality $< \mu\}$ ordered by inclusion, let $\Bbb Q =
\Bbb Q_1 \times \Bbb Q_2$.  Let $G = G_1 \times G_2 \subseteq \Bbb Q$ be
generic over $\bold V$ and let 
$\bold V_0 = \bold V,\bold V_1 = \bold V[G_1]$ and
$\bold V_2 = \bold V[G] = \bold V_1[G_2]$.
\mr
\item "{$\boxtimes_0$}"  In $\bold V_2,\underset\tilde {}\to B[G_1]$
is a Boolean Algebra of cardinality $\lambda$ with irr$^+(B) = \lambda$
and notational simplicity with a set of elements $\lambda$. \nl
[Why?  In $\bold V_1,\underset\tilde {}\to B[G_1]$ is like that by
\scite{1.3}.
Now as in $\bold V_1,\Bbb Q_2$ satisfies the $(\lambda,n)$-Knaster for every $n$
clearly by \scite{1.5} we are done.]
\ermn
In $\bold V_2$ the cardinal $\mu$ 
is still supercompact, hence it is well known that
\mr
\item "{$\boxtimes_1$}"  for every $Y \subseteq 2^\mu$ for some
normal ultrafilter ${\Cal D}$ on $\mu$ and $\bar Y = \langle Y_i:i < \mu
\rangle,Y_i \subseteq 2^{|i|}$ we have $\bar Y/{\Cal D}$ is $Y$ 
(i.e. $\bar Y/{\Cal D} \in \bold V^\mu_2/D$ and in the Mostowski Collapse of 
$\bold V^\mu_2/{\Cal D}$ the element $\bar Y/{\Cal D}$ is mapped to $Y$), 
hence $(2^\mu,Y,\mu,<)$ is isomorphic to $\dsize \prod_{i < \mu} 
(2^{|i|},Y_i,i,<)/{\Cal D}$.
\ermn
Again it is well known and follows from $\boxtimes_1$ that there is a 
sequence $\bar{\Cal D}^0 = \langle
{\Cal D}^0_\zeta:\zeta < (2^\mu))^+ \rangle$ of normal (fine) ultrafilters 
on $\mu$ satisfying: for each $\zeta < (2^\mu)^+$ the sequence 
$\bar{\Cal D}^0 \restriction \zeta$ belongs to (the
Mostowski collapse of) $\bold V^\mu_2/{\Cal D}_\zeta$.  In $\bold V_2$
we can code $\underset\tilde {}\to B = \underset\tilde {}\to B[G_1]$ and
${\Cal P}(\mu)$ and $\bar{\Cal D}^0 \restriction \kappa$ as a subset $Y$ of
$2^\mu = \lambda^*$ and get ${\Cal D},\bar Y$ as in $\boxtimes_1$ hence for some set
$A \in {\Cal D}$ of strongly inaccessible cardinals $> \chi$ 
there is a sequence 
$\langle (\mu_i,\lambda_i,B_i,\lambda^*_i):i \in A \rangle$ such that:
\mr
\item "{$(*)_1$}"  for $i \in A$ we have 
$i = \mu_i < \lambda_i \le \lambda^*_i = (\lambda^*_i)
^{\mu_i},\lambda_i$ is weakly inaccessible, Mahlo, $B_i$ is a Boolean 
Algebra generated by $\{x_\alpha:\alpha < \lambda_i\},
x_\alpha \notin \langle \{x_\beta:\beta < \alpha\} \rangle_{B_i}$,
irr$^+(B_i) = \lambda_i$ and, for notational simplicity, its sets of
elements is $\lambda_i$ 
\sn
\item "{$(*)_2$}"  $B$ is isomorphic to $\dsize \prod_{i \in A} 
B_i/{\Cal D}$ and $(\lambda^*,<) \cong \dsize \prod_{i \in A} 
(\lambda^*_i,<)/{\Cal D}$.
\ermn
Let $A^* = \{i < \mu:i$ strong inaccessible $> \chi\}$.  For $i \in \mu
\backslash A$ choose $\mu_i,\lambda_i,\lambda^*_i,B_i$ such that $(*)_1$
holds so $\mu_i = i$; why are there such $\lambda_i,B_i$?  Just e.g. use
$\lambda_{\text{Min}(A \backslash i)},B_{\text{Min}(A \backslash i)}$.
\sn
Let ${\Cal D}_i = {\Cal D}^0_i$ for $i < \kappa$ and ${\Cal D}_\kappa$ be
the ${\Cal D}$ as above.
So ${\Cal D}_i$ (for $i \le \kappa)$ is a normal ultrafilter on $\mu$ 
and we have
$i < j \le \kappa \Rightarrow {\Cal D}_i \in \bold V^\mu_2/{\Cal D}_j$, 
that is,
there is $\bar g = \langle g_{i,j}:i < j \le \kappa \rangle,g_{i,j} \in
{}^\mu({\Cal H}(\mu))$ such that ${\Cal D}_i$ is (the Mostowski collapse of)
$g_{i,j}/{\Cal D}_j \in \bold V^\mu_2/{\Cal D}_j$. \nl
All this was in $\bold V_2 = \bold V[G]$.
So we have $\Bbb Q$-names ${\underset\tilde {}\to {\bar g}} =
\langle {\underset\tilde {}\to g_{i,j}}:i < j \le \kappa \rangle,
{\underset\tilde {}\to {\bar{\Cal D}}} = \langle
{\underset\tilde {}\to {\Cal D}_i}:i \le \kappa \rangle$ and
$\langle ({\underset\tilde {}\to \mu_i},
{\underset\tilde {}\to \lambda_i},{\underset\tilde {}\to B_i},
{\underset\tilde {}\to \lambda^*_i}):i < \mu \rangle$.
As $\Bbb Q = \Bbb Q_1 \times \Bbb Q_2,\Bbb Q_2$ satisfies the $\mu^+$-c.c.
and $\Bbb Q_1$ is $\mu^+$-complete \wilog \, $\bar g$ is from $\bold V[G_1]$, so as we
could have forced first with some $\{f \in \Bbb Q_2:\text{Dom}(f)
\subseteq B\},B \in [\lambda^*]^{\le \mu}$; \wilog \, $\bar g$ and 
$\langle (\mu_i,\lambda_i,B_i,\lambda^*_i):i < \mu \rangle$ belong to $\bold V$.
Let $\Bbb P({\underset\tilde {}\to {\bar{\Cal D}}},
{\underset\tilde {}\to {\bar g}})$ be (the $\Bbb Q$-name of the) Magidor
forcing for $({\underset\tilde {}\to {\bar{\Cal D}}},
{\underset\tilde {}\to {\bar g}})$ (see \cite{Mg4}).  Let $\langle
{\underset\tilde {}\to \mu_i}:i < \kappa \rangle$ be the $\Bbb P
({\underset\tilde {}\to {\bar{\Cal D}}},{\underset\tilde {}\to {\bar g}})$-name of
the increasing continuous $\kappa$-sequence converging to $\mu$ which the
forcing adds and we can
restrict ourselves to the case $\mu_0 > \chi$.  Clearly clauses
$(\alpha)-(\zeta)$ in the conclusion hold for 
$\Bbb P = \Bbb Q * \Bbb P({\underset\tilde {}\to {\bar{\Cal D}}},
{\underset\tilde {}\to {\bar g}})$.  Now
\mr
\item "{$\boxtimes_2$}"  in $\bold V_2$, if $p \in \Bbb P(\bar{\Cal D},
\bar g)$ and $p \Vdash ``\underset\tilde {}\to f \in \dsize \prod_{i < \kappa}
\lambda_{{\underset\tilde {}\to \mu_i}}"$ \ub{then} 
there are $q$, an extension
of $p$ in $\Bbb P(\bar D,\bar g)$ and $f \in \dsize \prod_{i \in A^*}
\lambda_i$ such that \nl
$q \Vdash_{{\Bbb P}(\bar{\Cal D},\bar g)} `` \{i < \kappa:\underset\tilde {}\to f(i) =
f({\underset\tilde {}\to \mu_i})\} \in D^*"$.
\ermn
[Why?  By the properties of $\Bbb P(\bar{\Cal D},\bar g)$ there are a pure
extension $q_0$ of $p$ in $\Bbb P(\bar{\Cal D},\bar g)$ and sequence
$\langle u_i:i < \kappa \rangle$ such that above $q_0$ we
have: $\underset\tilde {}\to f(i)$ depends just on $\langle
{\underset\tilde {}\to \mu_j}:j \in u_i \cup \{i\} \rangle$ where $u_i
\subseteq i$ is finite.  As $D^*$ is a normal ultrafilter on $\kappa$, for some
$a^* \in D^*$ and finite $u \subseteq \kappa$ we have 
$i \in a^* \Rightarrow u_i = u$.
So there is $q$ such that $\Bbb P(\bar{\Cal D},\bar g) \models q_0 \le
q$ and $q \Vdash ``{\underset\tilde {}\to \mu_j} = \mu^*_j"$ for $j
\in u$, and so $f$ is well defined.]

Let $G_3 \subseteq \Bbb P(\bar{\Cal D},\bar g)$ be generic over $\bold
V_2$ and $\bold V_3 = \bold V_2[G_3]$ and let $\mu_i = {\underset\tilde {}\to \mu_i}
[G_3]$ so really $\langle \mu_i:i < \kappa \rangle$ is generic for
$\Bbb P(\bar{\Cal D},\bar g)$.  Now by $\boxtimes_2$ it follows that:
\mr
\item "{$\boxtimes_3$}"  in $\bold V_3 = \bold V_2[G_3]$ we have
$$
B \cong \dsize \prod_{i < \kappa} B_{\mu_i}/D^*.
$$
[Why?  In $\bold V_2$ there is an isomorphism $F$ from $B$ onto
$\dsize \prod_{i < \mu} B_i/{\Cal D} = \dsize \prod_{i \in A^*}
B_i/{\Cal D}_\kappa$, so let $F(x) = f_x/{\Cal D}_\kappa$ with $f_x
\in \dsize \prod_{i \in A^*} \lambda_i$ for $x \in B$, i.e. $x \in
\lambda$. \nl
In $\bold V_3$ let $f'_x \in \dsize \prod_{i < \kappa} \lambda_{\mu_i}$ be $f'_x(i) =
f_x(\mu_i)$ and we define a function $F'$ from $B$, i.e. from
$\lambda$ to $\dsize \prod_{i < \kappa} B_i/D^*$ by $F'(x) =
f'_x/D^*$.  Now $B \in {\Cal D} \Rightarrow \{i < \kappa:\mu_i \in B\}
= \kappa$ mod $J^{bd}_\kappa$ by the definition $\Bbb P(\bar{\Cal
D},\bar g)$, so as $F$ is one to one also $F'$ is, and $F'$ commute
with the Boolean operations as $F$ does; last $F'$ is onto by $\boxtimes_2$.]
\sn
\item "{$\boxtimes_4$}"  if $i < \kappa$ then ${\Cal
H}(\mu_{i+1})^{\bold V_3}$ is
the same as ${\Cal H}(\mu_{i+1})^{{\bold V}^{{\Bbb P}_i}_0}$, for some
$\mu^+_i$-centered forcing notion from ${\Cal H}(\mu_{i+1})$ 
(hence this forcing notion is $\lambda_{\mu_i}$-Knaster).
\ermn
[Why?  Note that ${\Cal H}(\mu_j)^{{\bold V}_2} = {\Cal
H}(\mu_j)^{\bold V_0}$ for $j \le \kappa$.
Also for each $i < \kappa$ in $\bold V_0$ there are 
${\Cal D}^i_j$, a normal ultrafilter on $\mu_i$
such that $(\bar{\Cal D}^i,\bar g^i) = 
(\langle {\Cal D}^i_j:j \le i \rangle,\langle 
g_{j_1,j_2} \restriction \mu_i:j_1 < j_2 \le i \rangle) \in \bold V$ 
is as above, i.e. $j_1 < j_2 \le i \Rightarrow {\Cal D}^i_{j_1} =
g_{j_1,j_2}/{\underset\tilde {}\to {\Cal D}^i_{j_2}} \in
\bold V^{\mu_i}/{\underset\tilde {}\to {\Cal D}^i_{j_2}},g^i_{j_1,j_2} \in
{}^{\mu_i}({\Cal H}(\mu_i))$ so $\Bbb P(\bar{\Cal D}^i,\bar g^i)$ is as in
\cite{Mg4}, and for some $G_{3,i} \subseteq \Bbb P(\langle {\Cal D}^i_j:
j \le i \rangle,\langle g_{j_1,j_2} \restriction \mu_i:j_1 \le \mu_2 \le i
\rangle)$ generic over $\bold V_0$ (equivalently over $\bold V_2$) we have
$G_{3,i} \in \bold V_3$ and ${\Cal H}
(\mu_{i+1})^{{\bold V}_3} = {\Cal H}(\mu_{i+1})^{{\bold V}_2[G_{3,i}]} =
{\Cal H}(\mu_{i+1})^{{\bold V}_0[G_{3,i}]}$.  See
\cite{Mg4}.  As $\Bbb P(\bar D^i,\bar g^i)$ is $\mu_i$-centered, clearly
$\boxtimes_4$ follows.] 
\nl
So obviously (by \scite{1.5})
\mr
\item "{$\boxtimes_5$}"  in $\bold V_3$, for each $i < \kappa$ we have
$B_i$ is a Boolean Algebra of cardinality $\lambda_{\mu_i}$, irr$^+(B_{\mu_i}) = 
\lambda_{\mu_i},\lambda_{\mu_i}$ is weakly Mahlo.
\ermn
Also in $\bold V[G_1]$, the forcing notion $\Bbb Q_2$ satisfies the
$\lambda$-Knaster condition and in $\bold V_2 = \bold V[G_1,G_2]$, the forcing
notion $\Bbb P(\bar{\Cal D},\bar g)$ from \cite{Mg1} is 
$\mu$-centered hence satisfies the $\lambda$-Knaster hence
\mr
\item "{$\boxtimes_6$}"  in $\bold V_3,B$ is a Boolean Algebra of
cardinality $\lambda$, a Mahlo cardinal and irr$^+(B) = \lambda$.
\ermn
Now let $\Bbb R = \text{Levy}(\mu^+,< \lambda)^{\bold V} = \{f \in
\bold V:\text{Dom}(f)
\subseteq \{(\alpha,\gamma):\alpha < \lambda,\gamma < \mu^+\},|\text{Dom}(f)|
\le \mu$ and for $\gamma < \alpha$, we have $f(\alpha,\gamma) < 1 + \alpha\}$,
ordered by inclusion.  
Clearly $\Bbb R$ satisfies the $\lambda$-Knaster condition, is 
$\mu^+$-complete in $\bold V$ and also in $\bold V_1$.
Let $G_{\Bbb R}$ be generic over $\bold V_1$.  Now in $\bold V[G_1,
G_{\Bbb R}]$, the forcing notion $\Bbb Q_2$ has the same definition and same
properties.  Also
(as in \cite{MgSh:433}, \cite{ShSi:677}), in $\bold V[G_1,G_2,G_{\Bbb R}]$ the
${\Cal D}_i(i \le \kappa)$ are still normal ultrafilters on $\mu$ and the
definition of $\Bbb P(\bar{\Cal D},\bar g)$ gives the same forcing notion
with the same properties and add the same family of subsets to $\kappa$
(as ${\Cal P}(\kappa)^{\bold V[G_1,G_2]} =
{\Cal P}(\kappa)^{\bold V[G_1,G_2,G_{\Bbb R}]}$).

So $G_{\Bbb R}$ is a subset of $\Bbb R$ generic over $\bold V[G_1,G_2,G_3]$.
Also in $\bold V[G_1,G_2],\Bbb R$ satisfies the $\lambda$-Knaster condition and in
$\bold V[G_1,G_2,G_{\Bbb R}],\Bbb P(\bar{\Cal D},\bar g)$ is $\mu$-centered
hence satisfies the $\lambda$-Knaster condition.  
Let $\bold V_4 = \bold V_3[G_{\Bbb R}]$, so in
$\bold V_4$ all the conclusions above holds but $\lambda = \mu^{++}$ hence
irr$(B) = \mu^+$ whereas irr$^+(B)$ remains $\lambda = \mu^{++}$.  
So we are done.
\hfill$\square_{\scite{1.7}}$\margincite{1.7}
\enddemo
\bigskip

\proclaim{\stag{1.8} Claim}  1) In the theorem \scite{1.7} we can replace

\block  ``a Boolean Algebra $B$ of cardinality $\lambda$, irr$^+(B) = \lambda$"
by e.g. ``a $\lambda$-Souslin tree" \endblock

The ``$\lambda$ strongly inaccessible Mahlo" is needed just for applying
\scite{1.3}, etc, but for $\dsize \prod_{i < \kappa} B_i/D^* \cong B$
is not needed (any model $M$, with universe $\subseteq \lambda$ is O.K.) \nl
2) We can apply the proof above to the proof in \cite{Sh:128} hence to theorem in
logics with Magidor Malitz quantifiers.
\endproclaim
\bigskip

\demo{Proof}  Similar to \scite{1.7}.
\enddemo
\newpage

\head {\S2 Consistency for small cardinals} \endhead  \resetall \sectno=2
 % \resetall 
\bn 
Theorem \scite{2.1} generalizes \scite{1.7} in some ways.  First $D^*$,
instead of being a normal ultrafilter on $\kappa$ is just a normal filter
which is large in appropriate sense so later can be applied to the case
$\kappa = \aleph_1$ (after a suitable preliminary forcing).  Second, we deal
with a general model and properties.  Thirdly, the forcing makes $\mu$ to
$\beth_\kappa$ (and more)
\proclaim{\stag{2.1} Theorem}  Suppose
\mr
\item "{$(a)$}"  $\bold V$ satisfies GCH above $\mu$ (for simplicity)
\sn
\item "{$(b)$}"  $\kappa$ is regular uncountable, $\aleph_0 \le \theta \le
\kappa < \chi < \mu < \vartheta < \lambda \le \lambda^* =
(\lambda^*)^\mu$, say $\vartheta
= \mu^+$
\sn
\item "{$(c)$}"  $\mu$ is supercompact, Laver indestructible or just
indestructible $\lambda^*$-hypermeasure (see exactly \cite{GM})
\sn
\item "{$(d)$}"  $D^*$ is a filter on $\kappa$ including the clubs and if
$f$ is a pressing down function on $\kappa$ then for some $u \in [\kappa]
^{< \theta}$ we have $\{\delta < \kappa:f(\delta) \in u\} \in D^*$
\sn
\item "{$(e)$}"  $\Bbb Q_1$ is a $(< \mu)$-directed complete forcing,
$|\Bbb Q_1| \le \lambda^*$ and $\Vdash_{\Bbb Q_1} ``\underset\tilde {}\to M$ is a 
model with universe $\lambda$ and 
vocabulary $\underset\tilde {}\to \tau \in {\Cal H}(\chi)"$
\sn
\item "{$(f)$}"  $\Bbb R$ is a $\mu^{++}$-complete forcing notion of
cardinality $\le \lambda^*$
\sn
\item "{$(g)$}"  $\Bbb Q_2$ is the forcing of adding $\lambda^*$ \,
$\mu$-Cohen subsets to $\mu$ and $\Bbb Q = \Bbb Q_1 \times \Bbb Q_2$ \nl
(see below Definition \scite{2.1}(a)).
\ermn
\ub{Then} for some forcing notion $\Bbb P$ we have $\Bbb Q_1 \times 
\Bbb Q_2 \times \Bbb R \lessdot \Bbb P$ and in $\bold V^{\Bbb P}$:
\mr
\item "{$(\alpha)$}"  forcing with $\Bbb P$ collapse no cardinal except those
collapsed by $\Bbb Q_1 \times \Bbb R$, in fact $\Bbb P/(\Bbb Q_1 \times
\Bbb Q_2 \times \Bbb R)$ is $\vartheta^-$-centered
\sn
\item "{$(\beta)$}"  forcing with $\Bbb P$ add no subset of $\chi$, forcing
with $\Bbb P/\Bbb Q_1 \times \Bbb Q_2 \times \Bbb R$ satisfies
$\boxtimes^1_{\gamma,\mu,\vartheta,\lambda,\lambda^*}$ from Definition
\scite{2.1a} below as witnessed by $\langle {\underset\tilde {}\to \mu_i}:i <
\kappa \rangle$
\sn
\item "{$(\gamma)$}"  ${\underset\tilde {}\to \mu_i} =
{\underset\tilde {}\to \mu_i}[G_{\Bbb P}],\mu$ is strong limit of cofinality
$\kappa$ and $\langle \mu_i:i < \kappa \rangle$ is an increasing continuous
sequence of strong limit singulars with limit $\mu$ (and ${\Cal H}
(\mu_{i+1})$ satisfies a parallel of the statement $\boxtimes_4$ from the
proof of \scite{1.7}),
\sn
\item "{$(\delta)$}"  for each $i < \kappa$ we have $\mu_i < \lambda_i 
\le \lambda^*_i
= (\lambda^*_i)^{\mu_i}$ and $\mu_\kappa = \mu,\lambda_\kappa = \lambda,
\lambda^*_\kappa = \lambda^*$ and $(\mu_i,\lambda_i,\lambda^*_i)$ is quite
similar to $(\mu,\lambda,\lambda^*)$ (see proof), more specifically: in some
intermediate universe $\bold V_1$, for some normal ultrafilter ${\Cal D}$
on $\mu$ and $F,F_*:\mu \rightarrow \mu$ we have $\dsize \prod_{i < \mu}
(F(i),<)/D \cong (\lambda,<),\lambda_i = F(\mu_i)$ and
$\dsize \prod_{i < \mu}(F_*(i),<)/{\Cal D} \cong (\lambda^*,<)$ and 
$F_*(\mu_i) = \lambda^*_i$
and we have $\bar M = \langle M_i:i < \mu \rangle$ and
$M_i$ a model with universe $\lambda_i$ and vocabulary $\tau$; and 
$\dsize \prod_{i < \mu} M_i/{\Cal D} \cong M$
\sn
\item "{$(\varepsilon)$}"  for $i < \kappa$ we have $2^{\mu_i} = \lambda^*_i$
and $2^{\lambda^*_i} = \mu_{i+1}$
\sn
\item "{$(\zeta)$}"  $\dsize \prod_{i < \kappa} M_{\mu_i}/D^*$ is isomorphic
to $M$ if $D^*$ is a normal ultrafilter, in fact, \nl
$\{\langle f(\mu_i):i < \kappa \rangle/_{D^*}:f \in \bold V_1$ and
$f \in \dsize \prod_{i < \mu} F(i)\}$ is the universe of
$\dsize \prod_{i < \kappa} M_{\mu_i}/D^*$
\sn
\item "{$(\eta)$}"  for every $f \in \dsize \prod_{i < \kappa} M_i/D^*$ 
we can in $\bold V_1$ find $\varepsilon(f) < \theta$ and
$g_{f,\varepsilon} \in \dsize \prod_{i < \mu} F(i)$ for $\varepsilon <
\varepsilon(f)$ such that $\{i < \kappa:\dsize \bigvee_{\varepsilon <
\varepsilon(f)} f(i) = g_{f,\varepsilon}(\mu_i)\} \in D^*$
\sn
\item "{$(\theta)$}"  $\dsize \prod_{i < \kappa} (\lambda_i,<)/D^*$ is
$\lambda$-like linear ordering (not necessarily well ordering as possibly
$\theta > \aleph_0$)
\sn
\item "{$(\iota)$}"  if 
$D^*$ is a normal ultrafilter, $Q_1 = Q^1_{\mu,\lambda}$ (of \scite{1.1})
and $\Bbb R = { \text{\rm Levy\/}}(\mu,< \lambda)$, \ub{then} the 
conclusion on irr in \scite{1.7} holds.
\endroster
\endproclaim
\bigskip

\definition{\stag{2.1a} Definition}  1) We say
$\boxtimes_{\gamma,\mu,\vartheta,\lambda^*}(\Bbb Q)$ or we say
$\Bbb Q$ satisfies $\boxtimes_{\gamma,\mu,\vartheta,\lambda^*}$ (as
witnessed by $({\underset\tilde {}\to {\bar \mu}},{\Cal D})$ \ub{if}:
\mr
\widestnumber\item{$(viii)$}
\item "{$(i)$}"  $\Bbb Q$ is a forcing notion of cardinality $\le \lambda^*$
\sn
\item "{$(ii)$}"  $\Bbb Q$ satisfies the $\vartheta$-c.c.
\sn
\item "{$(iii)$}"  $\Bbb Q$ (i.e. forcing with $\Bbb Q$) add a sequence
$\langle {\underset\tilde {}\to \mu_i}:i < \gamma \rangle$ of cardinals
$< \mu$, strongly inaccessible in $\bold V$, strong limit in $V^{\Bbb Q}$
\sn
\item "{$(iv)$}"  $\Vdash_{\Bbb Q} ``{\underset\tilde {}\to \mu_i} \,
(i < \gamma)$ is increasing continuous"
\sn
\item "{$(v)$}"  ${\Cal D}$ is a normal ultrafilter on $\mu$
\sn
\item "{$(vi)$}"  for every $p \in \Bbb Q$ for some $\beta < \gamma$ for
$A \in {\Cal D}$ there is $q$ satisfying $p \le q \in \Bbb Q$ such that 
$q \Vdash ``\{\mu_i:\beta < i < \gamma\} \subseteq A"$
\sn
\item "{$(vii)$}"  if $\gamma$ is a limit ordinal then
$\Vdash_{\Bbb Q} ``\mu = \dbcu_{i < \gamma} {\underset\tilde {}\to \mu_i}"$
\sn
\item "{$(viii)$}"  in $\bold V^{\Bbb Q}$ we have $2^\mu = \lambda^*$ and $\mu$ is
strong limit.
\ermn
2) We say $\boxtimes^+_{\gamma,\mu,\vartheta,\lambda^*}(\Bbb Q)$ or 
we say $\Bbb Q$ satisfies
$\boxtimes^+_{\gamma,\mu,\vartheta,\lambda^*}$ 
(as witnessed by $({\underset\tilde {}\to {\bar \mu}},f_\theta,
f_{\lambda^*})$ \ub{if}:
\mr
\item "{$(a)$}"  $\Bbb Q$ satisfies
$\boxtimes_{\gamma,\mu,\vartheta,\lambda^*}$ as witnessed by 
${\underset\tilde {}\to {\bar \mu}} = \langle
{\underset\tilde {}\to \mu_i}:i < \gamma \rangle$
\sn
\item "{$(b)$}"  if $G \subseteq \Bbb Q$ is generic over $\bold V$ \ub{then} for
every $\beta < \gamma$ we have ${\Cal H}(\mu_{\beta +1})^{\bold V^{\Bbb Q}}$ is
gotten from ${\Cal H}(\mu_{\beta +1})^{\bold V}$ by a forcing $\Bbb Q_{\beta +1}$
which is like $\Bbb Q$ for $\beta$.
\endroster
\enddefinition
\bigskip

\demo{Proof}  Like the proof of \scite{1.7} but we use \cite{GM} instead
of \cite{Mg4}; note that $\vartheta = \mu^{+3}$ comes from making the
forcing $\mu^{+3}$-c.c.  So the pure decision of $\Bbb P(\bar{\Cal D},\bar g)$
is changed accordingly.  Of course, the change in the assumption on
$D^*$ also has some influence.  \nl
${{}}$  \hfill$\square_{\scite{2.1}}$\margincite{2.1}
\enddemo
\bn
So we get e.g. \nl
\ub{\stag{2.2} Conclusion}:  Assume $\bold V$ satisfies ZFC + $\mu$ is a 
supercompact $+ ``\lambda > \mu$ is strong inaccessible".
\nl
1) For some forcing extension $\bold V^*$, for some ultrafilter $D^*$ on 
$\omega_1$ there is $\langle \lambda_i:i < \omega_1 \rangle$ such that:
\mr
\widestnumber\item{$(iii)$}
\item "{$(i)$}"   for $i < \omega_1,\lambda_i$ is
weakly inaccessible $< \beth_{\omega_1}$
\sn
\item "{$(ii)$}"  $\lambda = \beth^{++}_{\omega_1}$
\sn
\item "{$(iii)$}"  the linear order $\dsize \prod_{i < \omega_1} 
(\lambda_i,<)/D^*$ is $\lambda$-like, $\lambda_i$
first weakly inaccessible $> \beth_i$ (or first Mahlo $> \beth_i$).
\ermn
2) In part (1) we have: for some sequence 
$\langle B_i:i < \omega_1 \rangle$ of Boolean Algebras, each of cardinality
$< \beth_{\omega_1}$ we have Length$(\dsize \prod_{i < \omega_1} B_i/D^*) <
\dsize \prod_{i < \omega_1} \text{ Length}(B_i)/D^*$. \nl
3) If $\lambda$ in $\bold V,\lambda > \mu$ is Mahlo, also with irr.
\bigskip

\demo{Proof}  1) We start getting by forcing using a forcing notion from ${\Cal H}(\mu)$ 
(see \cite[Ch.XVI,2.5,p.793]{Sh:f} and history there) a normal filter 
$D^0$ on $\omega_1$ such that ${\Cal P}
(\omega_1)/D^*$ is layered \footnote{it means that this Boolean
Algebra is $\dbcu_{i < \omega_2} B^*_i,B^*_i$ is a Boolean Algebra of
cardinality $\aleph_1$, increasing continuous with $i$, and cf$(i) =
\aleph_1 \Rightarrow B_i \lessdot {\Cal P}(\omega_1)/D^*$}
and $\diamondsuit_{\aleph_1} + 2^{\aleph_1} =
\aleph_2$.  Hence (see \cite{FMSh:252} and history there) there is an 
ultrafilter $D^*$ on $\omega_1$ extending $D$ as required in \scite{2.1} 
clause (d) for $\theta = \aleph_1$, 
that is: if $g \in {}^{\omega_1}\omega_1$ is pressing down on
some member of $D$ then for some $\alpha < \omega_1,\{\beta < \omega_1:
g(\beta) < \alpha\} \in D$.  Now apply \scite{2.1} with $\theta = \aleph_1,
\Bbb R = \text{ Levy}(\mu^+,< \lambda),\lambda$ inaccessible. \nl
2) The proofs in \cite{MgSh:433} applies also in our changed
circumstances. \nl
3)  But for irr the problem seems more involved.  We use \scite{2.3} below
instead of \scite{1.3} and note that $\Bbb Q_2,\Bbb R$ and the Gitik
Magidor forcing $\Bbb P/(\Bbb Q_1 \times \Bbb Q_2 \times 
{\underset\tilde {}\to {\Bbb R}})$ though not fully preserving
$(*)_{\lambda,< \mu,{\underset\tilde {}\to B}}$ of \scite{2.3} below it
still leaves preserved for us 
$(*)_{\lambda,\aleph_0,{\underset\tilde {}\to B}}$ which is enough as we
now prove.  So in $\bold V^{\Bbb P}$ let $f_\alpha/D^* \in
\dsize \prod_{i < \kappa} B_{\mu_i}/D^*$ so
$f_\alpha \in \dsize \prod_{i < \kappa} B_{\mu_i}$ for $\alpha < \lambda$.  
For each $\alpha$ we can find in
$\bold V_2$ a sequence $\langle g_{\alpha,n}:n < \omega \rangle,
g_{\alpha,n} \in \dsize \prod_{i < \mu} B_i$ such that 
$\{i < \omega_1:(\exists n)(f_\alpha
(i) = g_{\alpha,n}(\mu_i))\} \in D^*$.  Without loss of generality we have
$A_{\alpha,n} = A_n$ where $A_{\alpha,n} = \{i < \omega_1:f_\alpha(i) =
g_{\alpha,n}(\mu_i)\}$, as $2^{\aleph_1} < \beth_{\omega_1} < \lambda =
\text{ cf}(\lambda)$.  Now in $\bold V_1$, there is an isomorphism $\bold j$ 
from $\dsize \prod_{i < \mu} B_i/{\Cal D}$ onto $B$, so $\bold j(g_{\alpha,n}/
{\Cal D}) \in B$.  In $\bold V_2[G_{\Bbb R}]$ we apply 
$(*)_{\lambda,\aleph_0,B}$
and find $\beta_0 < \beta_1 < \beta_2 < \beta_3 < \lambda$ 
such that $n < \omega \Rightarrow \bold j(g_{\beta_0,n}/{\Cal D}) = 
\sigma(\bold j (g_{\beta_0,n}/{\Cal D}),\bold j
(g_{\beta_0,n}/{\Cal D}),\bold j(g_{\beta_3,n}/{\Cal D}))$ where $\sigma$ is the
Boolean term $\sigma^*(x_0,x_1,x_2) = (x_0 \cap x_1) \cup (x_0 \cap x_2) \cup
(x_1 \cap x_2)$.  Hence

$$
Y_n =^{df} \{\zeta < \mu:B_\zeta \models g_{\beta_0,n}(\zeta) = \sigma^*
(g_{\beta_1,n}(\zeta),g_{\beta_2,n}(\zeta),g_{\beta_3,n}(\zeta))\} \in
{\Cal D}
$$
\mn
hence $Y = \dbca_{n < \omega} Y_n \in {\Cal D}$ hence for some $i^* < \kappa,
(\forall i)[i^* \le i < \kappa \rightarrow \mu_i \in Y]$ but $\mu_i \in Y
\Rightarrow (\forall n < \omega)[B_{\mu_i} \models 
g_{\beta_0,n}(\mu_i) = \sigma(g_{\beta_1,n}(\zeta),
g_{\beta_i,n}(\zeta),g_{\beta_3,n}(\zeta))]$.  As
$A_{\beta_\ell,n} = A_n$ we are done.  \hfill$\square_{\scite{2.2}}$\margincite{2.2}
\enddemo 
\bigskip

\remark{\stag{2.2a} Remark}  1) In 
\scite{2.2}(1),(2) \wilog \,\, $\beth_{\omega_1}$ is the
limit of the first $\omega_1$ (weakly) inaccessible. \nl
2) In \scite{2.2}(3) \wilog \, $\beth_{\omega_1}$ is the limit of the first 
$\omega_1$ Mahlo (weakly) inaccessible.  Can we omit Mahlo? \nl
3) Of course, \scite{2.2} is just one extreme variant.
\endremark
\bigskip

\proclaim{\stag{2.3} Claim}  1) For $\Bbb Q = \Bbb Q^1_{\mu,\lambda},
\underset\tilde {}\to B$ as in \scite{1.3} we have, for $\tau < \mu$ it is
forced $(\Vdash_{{\Bbb Q}^1_{\mu,\lambda}})$ that:
\mr
\item "{$(*)_{\lambda,\tau,{\underset\tilde {}\to B}}$}"  if
$y_{\alpha,\varepsilon} \in \underset\tilde {}\to B$ for $\alpha < \lambda,
\varepsilon < \tau$ \ub{then} for some $\beta_0 < \beta_1 < \beta_2 < \beta_3$
we have \nl
$\varepsilon < \tau \Rightarrow y_{\beta_0,\varepsilon} = \sigma^*
(y_{\beta_1,\varepsilon},y_{\beta_2,\varepsilon},y_{\beta_3,\varepsilon})$
where $\sigma^*(y_1,y_2,y_3) = (y_1 \cap y_2) \cup (y_1 \cap y_3) \cup (y_2
\cap y_3)$.
\ermn
2) If $B$ is a Boolean Algebra, $\tau < \lambda$ and $\Bbb Q^*$ is
$\tau^+$-complete (or just do not add new $\tau$-sequence of ordinals
$< |B|$) and satisfies the $(\lambda,4)$-Knaster property (i.e. among any
$\lambda$ conditions there are $\lambda$, any three of them has a common
upper bound), \ub{then} forcing by $\Bbb Q^*$ preserve
$(*)_{\lambda,\tau,B}$.   
\endproclaim
\bigskip

\demo{Proof}  1) As in \scite{1.3}
again the point is checking $(*)_{\lambda,\mu,{\underset\tilde {}\to B}}$ so let 
$p \Vdash ``\langle {\underset\tilde {}\to y_{\beta,\varepsilon}}:\beta <
\lambda,\varepsilon < \tau \rangle"$ be a counterexample.  For each $\alpha <
\lambda$ choose $p_\alpha$ such that $p \le p_\alpha$ and $p_\alpha \Vdash
``{\underset\tilde {}\to y_{\alpha,\varepsilon}} = y_{\alpha,\varepsilon}"$ 
for $\varepsilon < \tau$ and without loss of generality
$y_{\alpha,\varepsilon} \in p_\alpha$ and choose $\alpha_{\beta,\zeta}
\in w[p_\beta]$ 
for $\zeta < \mu$ such that $y_{\beta,\varepsilon} \in 
\langle \{x_\gamma:\gamma \in \{\alpha_{\beta,\varepsilon}:\varepsilon
< \zeta_\beta\} \rangle_{B[p_{\beta_\ell}]}$ for some $\zeta_\beta <
\tau^+$ with $\alpha_{\beta,\varepsilon}$ increasing with
$\varepsilon$, and let $\xi_\beta \le \zeta_\beta$ be such that
$(\forall \varepsilon)[\alpha_{\beta,\varepsilon} < \beta \equiv
\varepsilon < \xi_\beta]$.  Let $y_{\beta,\varepsilon} =
\sigma_{\beta,\varepsilon}(\ldots,x_{\alpha_{\beta,\varepsilon}},\ldots)_{\varepsilon
< \zeta_\beta}$ (so the term $\sigma_{\beta,\varepsilon}$ uses only
finitely many of its variables).  We choose
$S,w,r$, etc., as in the proof there with $\xi \le \zeta,\langle
\alpha_\varepsilon:\varepsilon < \xi \rangle,\langle
\sigma_\varepsilon:\varepsilon < \tau \rangle$ replacing $m \le
n,\langle \alpha_\ell:\ell < m \rangle,\sigma$.

We choose $\beta_0 < \beta_1 < \beta_2 < \beta_3$ in $S$ and it is enough
to find $q \in \Bbb Q^1_{\mu,\lambda}$ 
such that $\ell < 4 \Rightarrow p_{\beta_\ell} \le q$ and
$q \Vdash ``y_{\beta_0,\varepsilon} = \sigma(y_{\beta_1,\varepsilon},
y_{\beta_2,\varepsilon},y_{\beta_3,\varepsilon})$ for $\varepsilon < \tau"$.
We define $u = \dbcu_{\ell < 4} w[p_{\beta_\ell}]$ and ${\Cal F}$ as there,
i.e., 

$$
\align
\bigl\{ f:&f \in {}^u 2,f \restriction w[p_{\beta_\ell}] \in {\Cal F}
[p_{\beta_\ell}] \text{ for } \ell < 4 \text{ and for some} \\
  &\ell \in \{1,2,3\} \text{ we have} \\
  &m \in \{0,1,2,3\} \backslash \{\ell\} \and \zeta < \mu \Rightarrow
f(x_{\alpha_{\beta_0},\zeta}) = f(x_{\alpha_{\beta_m},\zeta}) \bigr\}.
\endalign
$$
\mn
Now check. \nl
2) Straightforward.   \hfill$\square_{\scite{2.3}}$\margincite{2.3}
\enddemo

%% you may want to move the following lines up a bit
\newpage
    
REFERENCES.  
\bibliographystyle{lit-plain}
\bibliography{lista,listb,listx,listf,liste}

\def\germ{\frak} \def\scr{\cal}
  \ifx\documentclass\undefinedcs\def\rm{\fam0\tenrm}\fi%f**k-amstex!
  \def\defaultdefine#1#2{\expandafter\ifx\csname#1\endcsname\relax
  \expandafter\def\csname#1\endcsname{#2}\fi} \defaultdefine{Bbb}{\bf}
  \defaultdefine{frak}{\bf} \defaultdefine{mathbb}{\bf}
  \defaultdefine{mathcal}{\cal}
  \defaultdefine{beth}{BETH}\defaultdefine{cal}{\bf} \def\bbfI{{\Bbb I}}
  \def\mbox{\hbox} \def\text{\hbox} \def\om{\omega} \def\Cal#1{{\bf #1}}
  \def\pcf{pcf} \defaultdefine{cf}{cf} \defaultdefine{reals}{{\Bbb R}}
  \defaultdefine{real}{{\Bbb R}} \def\restriction{{|}} \def\club{CLUB}
  \def\w{\omega} \def\exist{\exists} \def\se{{\germ se}} \def\bb{{\bf b}}
  \def\equivalence{\equiv} \let\lt< \let\gt> \def\cite#1{[#1]}
  \def\implies{\Rightarrow}
\begin{thebibliography}{FMSh 252}
\makeatletter \renewcommand{\@biblabel}[1]{[#1]} \makeatother
\par\noindent [References of the form {\tt math.XX/$\cdots$} refer to the
{\tt xxx.lanl.gov} archive]  \par

\bibitem[FMSh 252]{FMSh:252}Matthew Foreman, Menachem Magidor, and Saharon
  Shelah.
\newblock {Martin's maximum, saturated ideals and nonregular ultrafilters. II}.
\newblock {\em {Annals of Mathematics. Second Series}}, {\bf 127}:521--545,
  1988.

\bibitem[GM]{GM}Moti Gitik and Menachem Magidor.
\newblock {The Singular Cardinals Hypothesis revisited}.
\newblock In H.~Judah, W.~Just, and H.~Woodin, editors, {\em {The Proc. of MSRI
  conference on The Set Theory of the Continuum}}, Mathematical Sciences
  Research Institute Publications, pages 243--380. Springer Verlag, 1992.

\bibitem[Mg1]{Mg1}Menachem Magidor.
\newblock {On the singular cardinals problem II}.
\newblock {\em Annals Math.}, {\bf 106}:517--547, 1977.

\bibitem[Mg4]{Mg4}Menachem Magidor.
\newblock {Changing cofinality of cardinals}.
\newblock {\em Fund. Math.}, {\bf XCIX}:61--71, 1978.

\bibitem[MgSh 433]{MgSh:433}Menachem Magidor and Saharon Shelah.
\newblock {Length of Boolean algebras and ultraproducts}.
\newblock {\em {Mathematica Japonica}}, {\bf 48}(2):{301--307}, 1998.
 {\tt math.LO/9805145}

\bibitem[M2]{M2}J.~Donald Monk.
\newblock {\em {Cardinal Invariants of Boolean Algebras}}, volume 142 of {\em
  {Progress in Mathematics}}.
\newblock Birkh\"auser Verlag, Basel--Boston--Berlin, 1996.

\bibitem[Ru83]{Ru83}Matatyahu Rubin.
\newblock {A Boolean algebra with few subalgebras, interval Boolean algebras
  and retractiveness}.
\newblock {\em Trans. Amer. Math. Soc.}, {\bf 278}:65--89, 1983.

\bibitem[Sh:e]{Sh:e}Saharon Shelah.
\newblock {\em {Non--structure theory}}, accepted.
\newblock {Oxford University Press}.

\bibitem[Sh 128]{Sh:128}Saharon Shelah.
\newblock {Uncountable constructions for B.A., e.c. groups and Banach spaces}.
\newblock {\em {Israel Journal of Mathematics}}, {\bf 51}:273--297, 1985.

\bibitem[Sh:f]{Sh:f}Saharon Shelah.
\newblock {\em {Proper and improper forcing}}.
\newblock {Perspectives in Mathematical Logic}. {Springer}, 1998.

\bibitem[ShSi 677]{ShSi:677}Saharon Shelah and Otmar Spinas.
\newblock {On incomparability and related cardinal functions on ultraproducts
  of Boolean algebras}.
\newblock {\em Mathematica Japonica}, {\bf accepted}.
 {\tt math.LO/9903116}

\end{thebibliography}

\shlhetal

\enddocument%%

\bye